\documentclass[12pt]{article}  

\usepackage{hyperref}
\usepackage{amssymb}

\usepackage{amsfonts}
\usepackage{amsmath}

\usepackage{enumerate}

\newtheorem{theorem}{Theorem}

\newtheorem{proposition}[theorem]{Proposition}
\newtheorem{lemma}[theorem]{Lemma}
\newtheorem{corollary}[theorem]{Corollary}
\newtheorem{remark}[theorem]{Remark}

\newtheorem{definition}[theorem]{Definition}

\newcommand{\proof}{ \noindent{\it Proof:\ \ }}
\def\qed{\ifhmode\unskip\nobreak\fi\ifmmode\ifinner\else\hskip5 pt
\fi\fi\hbox{\hskip5 pt \vrule width4 pt height6 pt depth1.5 pt
\hskip 1pt }}
\newcommand{\po}{{\hspace*{-1ex}}{\bf .  }}

\let\oldmarginpar\marginpar
\renewcommand\marginpar[1]{${}^\clubsuit$\oldmarginpar[\raggedleft\scriptsize\sf #1]{\raggedright\scriptsize\sf #1}}
\usepackage{titlesec}
\titleformat{\section}{\normalfont\sffamily\Large}{\thesection.}{0.7em}{}
\titleformat{\subsection}{\normalfont\sffamily\filcenter}{}{0.7em}{}

\newcommand{\cref}[1]{Corollary~\ref{#1}}
\newcommand{\dref}[1]{Definition~\ref{#1}}

\newcommand{\lref}[1]{Lemma~\ref{#1}}
\newcommand{\pref}[1]{Proposition~\ref{#1}}
\newcommand{\rref}[1]{Remark~\ref{#1}}
\newcommand{\sref}[1]{Section~\ref{#1}}
\newcommand{\tref}[1]{Theorem~\ref{#1}}

\def\Bbb#1{\mathbb#1}
\newcommand{\R}{\Bbb{R}}
\newcommand{\V}{\Bbb{V}}

\newcommand{\N}{\Bbb{N}}
\newcommand{\Sp}{\Bbb{S}}

\newcommand{\C}{\Bbb{C\,}}

\newcommand{\isim}{isometric immersion }
\newcommand{\isims}{isometric immersions }
\newcommand{\na}{\nabla}

\newcommand{\ra}{\rangle}
\newcommand{\la}{\langle}

\newcommand{\spa}{{\rm span \,}}
\newcommand{\im}{{\rm Im \,}}
\newcommand{\rk}{{\rm rank \,}}

\newcommand{\Hess}{{\rm Hess}}

\newcommand{\D}{{\cal S}}

\newcommand{\moduli}{{\cal D}_f}
\newcommand{\modulii}{{\cal D}_i}
\newcommand{\CC}{{\cal C}}
\newcommand{\RR}{{\cal R}}

\begin{document}

\title{Classification of codimension two deformations\\
of rank two Riemannian manifolds}  
\author {Luis A. Florit\ \ \ and\ \ \ Guilherme M. de Freitas}

\date{}
\maketitle

\begin{abstract}  

The purpose of this work is to close the local deformation problem of
rank two Euclidean submanifolds in codimension two by describing their
moduli space of deformations. In the process, we provide an explicit
simple representation of these submanifolds, a result of independent
interest by its applications. We also determine which deformations are
genuine and honest, allowing us to find the first known examples of
honestly locally deformable rank two submanifolds in codimension two. In
addition, we study which of these submanifolds admit isometric immersions
as Euclidean hypersurfaces, a property that gives rise to several
applications to the Sbrana-Cartan theory of deformable Euclidean
hypersurfaces.

\end{abstract}

\section*{}  

Among the most fundamental properties of a structure defined on a certain
class of objects is its rigidity, in a broad sense: whether the structure
exists on the given class at all, when it does if it is unique, or, when
not unique, to somehow understand its moduli space, i.e., the space of
deformations of the structure. In submanifold theory, the corresponding
concept is that of isometric rigidity, and starts by asking if a fixed
Riemannian manifold admits, either locally or globally, an isometric
immersion in a given ambient space, usually the Euclidean space.
And, when it does exist, to try to classify all its \isims in that
ambient space.

Well-known results in the field include the Nash-Gromov-Rocklin embedding
theorems, which states that any smooth Riemannian manifold $M^n$ admits
an \isim in the Euclidean space $\R^{n+p}$ with codimension $p\sim
n^2/2$. This implies that the rigidity question for submanifolds only
makes sense for relatively small codimensions, at least bounded in terms
of the dimension of the manifold. One of the main characteristics of
rigidity problems in submanifold theory is that, except for the lowest
dimensions, the difficulties do not depend on the dimension of the
manifold, but rather they grow very fast with the codimension.

In this work we are interested in the local deformation problem of the
theory. A well-known result and the starting point of the subject is the
Beez-Killing theorem, that states that a Euclidean hypersurface is
locally rigid provided that the number of nonzero principal curvatures,
called the {\it rank} of the hypersurface, is at least three. This result
has had several generalizations, like the ones in \cite{alle},
\cite{bbg}, \cite{dd} and \cite{si}. The general idea behind these works
is that, in order for a submanifold in low codimension to possess
noncongruent deformations, its second fundamental form and curvature
tensor must be highly degenerate. Several of this kind of results were
unified and generalized in \cite{df2}, where the notion of {\it genuine
rigidity} was introduced (see Section 1 for definitions). This concept
relies on the idea that, as we discard congruent submanifolds when
analyzing rigidity, we should also discard deformations that are induced
by deformations of a bigger dimensional submanifold containing the
original one.

On the other hand, flat hypersurfaces have at most rank one, and the
moduli space of their local isometric deformations is well understood:
they are parametrized in a very simple and geometric way, through their
Gauss map, by smooth arbitrary regular curves in their ambient space
$\R^{n+1}$. In this work we will see that this is also the case for flat
submanifolds in codimension two. Therefore, we argue here that we should
also discard compositions with flat submanifolds since these are well
understood in low codimension, giving rise to the concept of {\it honest
rigidity}.

 From the discussion above we conclude that the interesting local
deformation phenomena for hypersurfaces arise only for those that have
rank two. A century ago, V. Sbrana and E. Cartan described the rank two
locally deformable hypersurfaces, the so-called {\it Sbrana-Cartan
hypersurfaces}, by showing that they split into four classes, according
to their space of deformations. However, despite their
classification, a deep question remained open for almost a century,
namely, the very existence of examples in the least deformable {\it
discrete} class, whose members admit, precisely, only one noncongruent
deformation. This question was answered for any dimension in
\cite{dft2} by means of a very geometric construction carried out in
codimension two that provided a large family of examples in this discrete
class. More precisely, it was shown that the transversal intersection of
two flat hypersurfaces in general position gives rise to a Sbrana-Cartan
hypersurface of the discrete class, together with its two unique
isometric immersions as a Euclidean hypersurface. We call these
hypersurfaces {\it Sbrana-Cartan hypersurfaces of intersection~type}.


\medskip

Since the problem beyond hypersurfaces is quite involved, until recently
nothing similar to the Sbrana-Cartan theory for codimension higher than
one had been attempted, not even in codimension two. For the compact
case, the main result in \cite{dm2} says that, if $n\geq 5$,
they are nowhere genuinely deformable once certain mild singularities for
the extensions are allowed. However, the necessity of considering
singular extensions was not established in \cite{dm2}, and we will
address this question here.

Now, in the search for interesting local deformation phenomena it is
natural to begin such study with the Euclidean submanifolds in
codimension two whose second fundamental forms or curvature tensors are
very degenerate, that is, those with rank two. We point out that the rank
condition, in this setting, is essentially an intrinsic property, and
agrees with the rank of the curvature tensor of the manifold.

When substantial and irreducible, rank two submanifolds in codimension
two naturally divide into three classes, {\it elliptic}, {\it parabolic}
and {\it hyperbolic}, according to the number 0, 1 or 2 of independent
normal directions whose shape operators have rank one. In particular, in
the hyperbolic case we have a fundamental function, the {\it main angle}
between these two normal directions, that plays a key role in this work.
To our surprise, elliptic and parabolic submanifolds were shown to be
honestly rigid in \cite{df3} without the need of a thorough understanding
of their space of deformations. However, the study of the deformations of
hyperbolic submanifolds remained elusive, among other reasons, by the
lack of a good representation of them.

\medskip

And this is precisely the main goal of this work: to describe locally in
a convenient way all these hyperbolic submanifolds (\tref{param}) and to
characterize their moduli spaces of deformations (\tref{mainthg}), hence
closing the classification of the local deformation problem of rank two
Euclidean submanifolds in codimension two. We will see that, in contrast
to the elliptic and parabolic cases, there are hyperbolic submanifolds
that are honestly deformable.

One of the main objectives of this work is to call the attention to the
other side of the coin of genuine and honest rigidity, that remained
somehow hidden: although it is a good idea to discard deformations that
arise from submanifolds in lower codimension when dealing with rigidity,
one can instead study non genuine deformations but in higher codimension
to obtain information about the deformations in the codimension we are
interested in. The last two sections of this paper show that this twist
is indeed fruitful. Therefore, we should not simply ignore all non honest
and non genuine deformations, but instead it is important to understand
them.

\medskip

To accomplish our task, we proceed as follows. First, in \sref{project}
we show that all the data relevant to the study of the deformations of a
hyperbolic submanifold in codimension two project to the leaf space of
its totally geodesic foliation of relative nullity, that is a smooth
surface. We then use this to show in \sref{psatpg} that these hyperbolic
submanifolds admit a very simple and geometric representation in terms of
their {\it polar surfaces}. This result is interesting in its own right,
beyond rigidity problems, since it provides an explicit tool to construct
interesting classes of submanifolds. For example, as an immediate
consequence of this parametrization, we obtain a local description of all
flat submanifolds in codimension two in terms of Euclidean surfaces with
flat normal bundle, which were classified in \cite{fera}. This
construction is not only simpler to use and more explicit, but also more
elegant than the one found in \cite{df0}.

In \sref{ihs} we provide a full set of 6 functions of a given hyperbolic
submanifold that are invariant under deformations, which will
substantially simplify our work. This is a somehow different approach
than the usual ones that deal with deformability of submanifolds. We then
make use of this set to obtain our main result in \sref{mmm}: the
description of the moduli space of deformations. It turns out that this
space is given as pairs of functions of one variable satisfying a single
equation that, surprisingly, depends only on the metric of the polar
surface. This is not to say that this space is always easy to compute.
Indeed, the richness of the deformation phenomena arises from the way
that these two functions are entangled in the equation.

We determine in \sref{gensec} which deformations are honest, and how the
ones that are not extend. As a simple consequence, we will conclude the
necessity of considering singular extensions for the global rigidity
problem for compact manifolds, as established in \cite{dm2}. In
\sref{hsahg} we characterize which hyperbolic submanifolds in codimension
two admit an \isim as a Euclidean hypersurface, a result that should be
seen as the converse to the main result in \cite{dft3}, and that also has
applications to the Sbrana-Cartan theory.

Finally, we use this machinery to give three applications. First, in
\sref{mangle} we characterize all deformations that preserve the main
angle. It turns out that all of them are not honest nor genuine, yet
they provide applications to the Sbrana-Cartan theory, in
particular, giving new examples of the interesting classes. Then, in
\sref{sbcait} we fully recover the main result in \cite{dft2} cited above
about the construction of the Sbrana-Cartan hypersurfaces of intersection
type. More importantly, we compute their moduli space of deformations as
submanifolds in codimension two, where they naturally live. We will see
that generically they are also honestly rigid, except in one situation
where the moduli space is precompact and diffeomorphic to a line. These
are the first known examples of honestly locally deformable rank two
Euclidean submanifolds in codimension two.

\medskip

We end this introduction by pointing out that the techniques used in this
work can be easily extended to elliptic and parabolic submanifolds in
codimension two by using complex conjugate coordinates instead of real
ones, in the same spirit as in \cite{dft3}. Therefore, a unified approach
of the rigidity and deformation phenomena of rank two submanifolds in
codimension two, and in fact in space forms, can be easily carried out
with the techniques presented here. However, since as we pointed out
elliptic and parabolic submanifolds are honestly rigid, by simplicity of
the presentation and to avoid conceptual duplications we choose to
restrict ourselves to the hyperbolic case in Euclidean space.

\section{Preliminaries}\label{prelg}  

In this section we recall basic facts about rank two submanifolds and
their splitting tensors, and we introduce the concept of honest rigidity,
which is a slight extension of that of genuine rigidity and better
suits the study of deformation phenomena when the submanifold admits
an \isim in lower codimension.

\medskip

Along this paper, $M^n$ will denote an $n$-dimensional Riemannian
manifold, and $f\colon M^n\to\R^{n+p}$ an isometric immersion, always in
codimension $p=2$ except in this section. By a {\it deformation} of $f$
we simply mean another non-congruent \isim of $M^n$ into the same ambient
space $\R^{n+p}$.

\subsection{Honest rigidity}  

To study the deformation problem for such an $f$ we need the following.
\begin{definition}\po\label{comp}
{\rm We say that $f:M^n\to\R^{n+p}$ is a {\it composition} (of $g$) if
there is an open subset $U\subset\R^{n+r}$ with $r<p$, and isometric
immersions $h:U\subset\R^{n+r}\to\R^{n+p}$ and $g:M^n\to \R^{n+r}$ with
$g(M^n)\subset U$ such that $f=h\circ g$.}
\end{definition}

It is important to observe that all flat Euclidean hypersurfaces, as $h$
above for $r=p-1$, can be easily parametrized and classified using just
an arbitrary curve in the unit sphere $\Sp^{n+1}\subset\R^{n+2}$ (the
Gauss image of $h$), together with a function of one variable, by means
of the Gauss parametrization; see e.g.\,(15) in \cite{dft2}. In fact, as
we will see in \cref{flat} below, in this paper we also describe easily
flat Euclidean submanifolds in codimension two. Therefore, compositions
are not interesting when studying rigidity in codimension two and three,
since all these deformations arise from the ones for the submanifold $g$,
hence reducing the codimension of the problem. This is one of the reasons
why we want to discard compositions, in the same way we discard
congruences.

\medskip

Key concepts in this work are those of genuine rigidity and genuine
deformations introduced in \cite{df2} and extended to the conformal realm
in \cite{ft}. As we pointed out, the complexity of rigidity problems in
submanifold theory grows very fast with the codimension. Moreover, if we
have inclusions $M^n\subset N^{n+k}\subset\R^{n+p}$, deformations of the
lower codimensional submanifold $N^{n+k}$ in $\R^{n+p}$ induce obvious
deformations of $M^n$ in $\R^{n+p}$. In order to discard these simpler
deformations we proceed as follows.

\medskip

Given another \isim $\hat f: M^n\to\R^{n+q}$, $q\in\N$, we say that
the pair $\{f.\hat f\}$ {\it extends isometrically\/} if there are an
isometric embedding $j:M^n\hookrightarrow N^{n+k}$ into a Riemannian
manifold $N^{n+k}$, $k\geq 1$, and a pair of \isims 
$F:N^{n+k}\to\R^{n+p}$ and $\hat F:N^{n+k}\to\R^{n+q}$,
such that $f=F\circ j$ and $\hat f=\hat F\circ j$,
that is, when the following diagram commutes: \vspace{1ex}

\begin{picture}(110,84)
\put(155,31){$M^n$}
\put(209,31){$N^{n+k}$}
\put(242,62){$\R^{n+p}$}
\put(242,0){$\R^{n+q}$}
\put(195,59){${}_f$}
\put(195,9){${}_{\hat f}$}
\put(240,46){${}_F$}
\put(239,26){${}_{\hat F}$}
\put(193,40.5){${}_j$}
\put(228,42){\vector(1,1){16}}
\put(228,28){\vector(1,-1){16}}
\put(177,44){\vector(3,1){60}}
\put(177,26){\vector(3,-1){60}}
\put(185,34){\vector(1,0){21}}
\put(185,36){\oval(7,4)[l]}
\put(205,43){$\circlearrowleft$}
\put(205,20){$\circlearrowright$}
\end{picture}
\vspace{1.8ex}

\noindent Accordingly, we say that $\hat f$ as above is a {\it genuine
deformation} of $f$ (in codimension $q$), or simply that $\{f,\hat f\}$
is a {\it genuine pair}, if $\{f|_U,\hat f|_U\}$ does not extend
isometrically on any open subset $U\subset M^n$.

\medskip

Observe that this concept (locally) extends that of compositions when
$k=p<q$, and that of congruence when $k=p=q$. This allowed to unify and
generalize in \cite{df2} and \cite{ft} several known rigidity results that
seemed different in nature. The main results in those two papers are
that, in sufficiently low codimensions, the members of a genuine pair
have to be mutually ruled, with a special kind of rulings of large
dimension.

\medskip

The associated rigidity concept is the following.

\begin{definition}\po\label{gr}
{\rm An \isim $f:M^n\to\R^{n+p}$ is {\it genuinely
rigid} if, for any given \isim $\hat f: M^n\to\R^{n+p}$,
there is an open dense subset $U\subset M^n$ such that the pair
$\{f|_U, \hat f|_U\}$ extends isometrically.}
\end{definition}

Although this concept is appropriate when $M^n$ admits no lower codimensional \isim $g$ as in
\dref{comp}, when it does, we automatically have all the compositions
$h\circ g$ for each flat submanifold $h$, and therefore we cannot expect
to have genuine rigidity. On the other hand, as we saw, these
compositions in low codimension are also well understood, and thus we
want to discard them as well. Therefore, this justifies us to say
that a genuine deformation $\hat f: M^n\to\R^{n+p}$ of $f$ is
{\it honest} if it is nowhere a composition.
Accordingly to this, we introduce our next concept.

\begin{definition}\po\label{hr}
{\rm We say that $f$ is {\it honestly rigid\/} if its only genuine
deformations are compositions along an open dense subset.}
\end{definition}

We do not require in this work for $M^n$ to admit no \isim as an
Euclidean hypersurface, as it was done in Theorem 1 in \cite{df3}. This
has two main reasons. First, we can and we will use the machinery
developed in this paper to obtain information about the theory of
deformable hypersurfaces; see Sections \ref{mangle} and \ref{sbcait} for
examples of this situation. Secondly, the concept of honest rigidity
better fits the deformation problem in codimensions bigger than one and
avoids this hypothesis. For example, the cited hypothesis in Theorem 1 in
\cite{df3} becomes now unnecessary, since its proof actually shows
the stronger result that,
indeed, any elliptic rank two Euclidean submanifold in codimension two is
honestly rigid, even if it is also a Euclidean hypersurface.
Analogously, Theorem 4 in \cite{df3} implies that any parabolic rank two
Euclidean submanifold in codimension two is honestly rigid. We will show
in \sref{sbcait} that this is not the case for hyperbolic submanifolds in
codimension two.

\subsection{Rank two and hyperbolic submanifolds}  

Let $R$ be the curvature tensor of $M^n$. We denote by $\Gamma(x)$ the
{\it nullity} of $M^n$ at $x\in M^n$,
$$
\Gamma(x):=\{X\in T_xM:R(X,Y)=0\mbox{ for all }Y\in T_xM\}.
$$
The {\it rank} of $M^n$ at $x$ is the integer $n-\dim\Gamma(x)$, which is
constant on connected components of an open dense subset of $M^n$. Now,
given an \isim $f:M^n\to\R^{n+2}$ of $M^n$, we denote by $\Delta(x)$ the
{\it relative nullity} of $f$ at $x$, that is, the nullity space of the
second fundamental form $\alpha=\alpha_f$ of $f$ at $x$,
$$
\Delta(x)=\{X\in T_xM:\alpha(X,Y)=0\mbox{ for all }Y\in T_xM\}.
$$
We call the  {\it rank} of $f$ at $x$ the integer $n-\dim\Delta(x)$.
By the Gauss equation for $f$ it is immediate that
$\Delta(x)\subset\Gamma(x)$, and so the rank of $f$ is pointwise greater
than or equal to the rank of $M^n$.
It is well-known that both $\Gamma$ and $\Delta$ are smooth, integrable
totally geodesic distributions on $M^n$ (along the connected components
of an open dense subset of $M^n$ where they have constant dimension).
In addition, $\Delta$ is totally geodesic also in the ambient Euclidean
space.

Since our work is local in nature, we will assume whenever necessary and
without further mention that all distributions that appear as images or
kernels of tensors have constant dimension. This will not bring any
problem since our rigidity concepts are required to hold locally almost
everywhere by definition.
In fact, we will see in the last section that submanifolds that deform in
very different forms can be glued smoothly, although non analytically, in
quite complicated ways, and then we cannot expect any local
classification to hold everywhere.

\medskip

In particular, the following proposition shown in \cite{df3} tells us
that we can restrict ourselves to rank two immersions when studying
deformations of rank two Riemannian manifolds.

\begin{proposition}\po\label{diffnul}
Let $f:M^n\to\R^{n+2}$ be an isometric immersion of a Riemannian manifold
of rank two. If at a certain point $x\in M^n$ it holds that
$\dim\Delta(x)\neq n-2$, then $f|_W$ is a composition in some open
neighborhood $W$ of $x$.
\end{proposition}

We denote by $N^1_f(x)$ the {\it first
normal space} of $f$ at $x$,
$$
N^1_f(x)=\spa\{\alpha(X,Y): X,Y\in T_xM\}\subset T^\perp_xM.
$$
If $M^n$ is nowhere flat, then by the Gauss equation $\dim N^1_f\geq 1$.
If, on the other hand, $\dim N^1_f = 1$, it is easy to see using Codazzi
equation that $N^1_f$ is parallel, and hence $f(M^n)$ is contained in an
affine hyperplane of $\R^{n+2}$, being, in particular, a composition.
But even if $M^n$ is flat, $\dim N^1_f < 2$ also implies that $f$ is a
composition of a totally geodesic inclusion $M^n\subset \R^{n+1}$ by
Proposition 9 in \cite{df2}. Therefore, we also have the following,
where $f$ to be {\it full} means that $N^1_f=T^\perp_fM$ everywhere.

\begin{proposition}\po\label{n12}
If $f:M^n\to\R^{n+2}$ is any isometric immersion that is nowhere a
composition, then $\dim N^1_f = 2$  almost everywhere.
\end{proposition}

As it is the case with any totally geodesic distribution, $\Delta$
possesses its \emph{splitting tensor}
$C:\Delta\times\Delta^\perp\rightarrow\Delta^\perp$ given by
$$
C(T,X)=C_TX=-(\nabla_XT)_{\Delta^\perp}.
$$
By the Codazzi equation, for any vector $\xi$ in the normal bundle of
$f$, the corresponding shape operator $A_\xi$ of $f$, that we always
consider restricted to $\Delta^\perp$, satisfies that
\begin{equation}\label{stcod}
A_\xi \circ C_S=C_S^{\,t}\circ A_\xi,
\end{equation}
for any $S\in\Delta$. In particular, when $M^n$ has rank two and
$f$ is nowhere a composition, by \pref{diffnul} we have $\Delta=\Gamma$,
the splitting tensor of $\Delta$ agrees with the one for the nullity
$\Gamma$ of $M^n$, and therefore it is intrinsic; see \cite{dft2} for
details.

\begin{remark}\po\label{rank3}
{\rm Given a rank two Riemannian manifold $M^n$ with splitting tensor $D$
of $\Gamma$, if $d:=\dim (\im D\subset {\rm End}(\Gamma^\perp))=4$,
condition \eqref{stcod} clearly implies that $M^n$ admits no rank~2
isometric immersion in Euclidean space. In particular, by \pref{diffnul},
$M^n$ cannot be a Euclidean submanifold in codimensions 1 or 2. On the
other hand, if $d=3$, also from~\eqref{stcod} we see that any rank~2
\isim $f$ of $M^n$ must satisfy $\dim N^1_f=1$, and thus by \pref{n12} it
must be a composition almost everywhere.}
\end{remark}

Recall that $f$ is {\it surface-like} if, along each connected component
$U_\lambda$ of an open dense subset $U$ of $M^n$, there is a surface
$V^2$ such that $U_\lambda$ splits as a Riemannian product
$U_\lambda\subset V^2\times\R^{n-2}$ if the splitting tensor $C$ vanishes
(respectively, $U_\lambda\subset CV^2\times\R^{n-3}$ if $C\neq 0$), and
$f|_{U_\lambda}=(g\times Id_{\R^{n-2}})|_{U_\lambda}$ splits
(respectively,
$f|_{U_\lambda}=({\cal C}g\times Id_{\R^{n-3}})|_{U_\lambda}$),
for some \isim $g:V^2\to\R^4$ (respectively, $g:V^2\to\Sp^4$, and
${\cal C}g$ stands for the cone over $g$, ${\cal C}g(x,t)=tg(x)$).
It is easy to see that $f$ is surface-like if and only if
it holds everywhere that $C_T=\mu(T)I$ for all
$T\in\Delta$; cf. Lemma 6 in \cite{dft2}.

\medskip

Now, for a rank two immersion $f$ as above, since its codimension is
two, for a given basis
$X,Y\in\Delta^\perp(x)$ there are $a,b,c\in\R$ such that
\begin{equation}\label{hpeg}
a\alpha_f(X,X)+2c\alpha_f(X,Y)+b\alpha_f(Y,Y)=0.
\end{equation}
Following \cite{df1}, a nowhere surface-like rank two Euclidean
submanifold $f:M^n\to\R^{n+2}$ is called \emph{hyperbolic} (respectively,
\emph{parabolic} or \emph{elliptic}) if it holds everywhere that
$\dim N^1_f = 2$, and $ab-c^2<0$ (respectively, $ab-c^2=0$ or
$ab-c^2>0$), a condition that is independent of the given basis.
Of course, these concepts make perfect sense in any codimension, but in
this paper we always reserve the term `hyperbolic (parabolic, elliptic)
submanifold' to those in codimension two, except for surfaces. We can
choose the basis $\{X_1,X_2\}$ such that \eqref{hpeg} takes the form
\begin{equation}\label{hpesg}
\alpha_f(X_1,X_1)-\epsilon\alpha_f(X_2,X_2)=0,
\end{equation}
where $\epsilon=1$ (respectively, $\epsilon=0$, $\epsilon=-1$).
Moreover, the pairs $a_1X_1+a_2X_2,\mbox{ }a_1X_2+\epsilon a_2X_1$ also
satisfy \eqref{hpesg} and, up to signs, there are no others.
Then let
$J:\Delta^\perp\rightarrow\Delta^\perp$ be the (unique up to sign) linear
map defined by $JX_1=X_2$ and $JX_2=\epsilon X_1$. In particular,
$J^2=\epsilon I$. We conclude that $f$ is hyperbolic (respectively,
parabolic or elliptic) if and only if there is a linear map
$J:\Delta^\perp\rightarrow\Delta^\perp$ such that $J^2=\epsilon I$ with
$\epsilon=1$ (respectively, $\epsilon=0$ or $\epsilon=-1$), and
$\alpha_f(X,JY)=\alpha_f(JX,Y)$, for all $X,Y\in\Delta^\perp$,
or, equivalently,
\begin{equation}\label{intcod}
A_\xi\circ J=J^{\,t}\circ A_\xi,\ \ \ \forall\ \xi\in T^\perp_fM.
\end{equation}
In particular, it is easy to check that $\epsilon+1$ is the number of
linearly independent normal directions whose corresponding shape
operators have rank one. Moreover, by \eqref{stcod},
\begin{equation}\label{stdhsg}
\{C_S:S\in\Delta\}\subset\mbox{span}\{I,J\}.
\end{equation}
Hence, since $f$ is always assumed to be nowhere surface-like, by
\eqref{stdhsg} the endomorphism~$J$ above is also intrinsic, and so is
the property of being hyperbolic, parabolic or elliptic. In other words,
by \rref{rank3}, if a nowhere surface-like rank two Riemannian
manifold~$M^n$ admits a Euclidean \isim in codimension two that is
nowhere a composition, then $d\leq 2$, $\Delta=\Gamma$, the immersion has
to be either hyperbolic, parabolic or elliptic, \eqref{stdhsg} holds, and
accordingly we call $M^n$ itself hyperbolic, parabolic or elliptic. This
justifies the intrinsic flavour of the title of this work. In particular,
if a deformation of a hyperbolic submanifold is not hyperbolic, then it
is somewhere a composition, and hence it is not honest.

\begin{remark}\po\label{gdh}
{\rm As we pointed out, it is very easy to classify all the compositions
in codimension two of a given hypersurface $M^n\subset\R^{n+1}$. On the
other hand, the genuine \isims of such an $M^n$ in codimension two
were classified in Proposition 3 and Theorem 4 of \cite{df2} for the
parabolic case, and in Theorem 1 of \cite{dft3} for the remaining
elliptic and hyperbolic cases.}
\end{remark}

\begin{remark}\po\label{z}
{\rm The local understanding of rank two Riemannian manifolds is also
important from an intrinsic point of view. Indeed, by Theorem A
in \cite{fz}, such a complete Riemannian manifold with finite volume is
always surface-like. This is true even allowing rank less or equal than
two. In particular, this shows that Nomizu's conjecture, which states
that a complete locally irreducible semi-symmetric space of dimension
at least three must be locally symmetric, that is well-known to be
false, is actually true for complete manifolds with finite volume.}
\end{remark}

\subsection{The Sbrana-Cartan theory}  


As we pointed out, a Euclidean hypersurface is locally rigid when its
rank is greater than or equal to three, and highly deformable but well
understood when flat, i.e., when its rank is at most one. At the
beginning of the last century, V. Sbrana \cite{sb} and a few years later
E. Cartan \cite{ca}, independently and with different techniques, worked
out the remaining interesting case, the hypersurfaces with rank two. That
is, they classified nowhere flat locally deformable Euclidean
hypersurfaces, $f:M^n\to\R^{n+1}$, for $n\geq 3$, extending earlier works
by Schur and Bianchi; see \cite{dft2} and references therein. According
to this classification, these hypersurfaces, now called
{\it Sbrana-Cartan hypersurfaces}, are (locally) divided into four
classes.

The first two classes of Sbrana-Cartan hypersurfaces are obvious and
highly deformable: the surface-like and the ruled ones. The first ones
deform as their surfaces do in $\R^3$ or $\Sp^3$, while the space of
deformations of a ruled hypersurface can be naturally parametrized by the
set of smooth real functions in one variable.

Therefore the actually interesting Sbrana-Cartan hypersurfaces belong
to the two remaining {\it continuous} and {\it discrete} classes,
and thus are the ones that demand the hard work in the theory.
The ones in the continuous class admit precisely a one parameter family
of deformations, while the ones in the discrete class have just
only one noncongruent deformation.

Since the beginning of the theory several families of examples of the
continuous class were known. In fact, the bulk of Sbrana and Cartan works
is concentrated on the study of this class. For example, those whose
Gauss map is a minimal surface in the sphere belong to this class.
Another large set of examples is given by minimal hypersurfaces of rank
two that have a one parameter associated family of deformations like
minimal surfaces do (\cite{dm1}).

However, until very recently not a single example of the discrete class
was known, nor even if this class was actually empty. A large set of
explicit examples of the discrete class was then explicitly constructed
and characterized in \cite{dft2} in a very geometric way: as the
transversal intersection of two generic flat hypersurfaces. Although the
construction is quite natural, the actual computations are long and
involved. It turns out that these submanifolds in codimension two are
hyperbolic, and in \sref{sbcait} we will recover this result, in a much
simpler way, by using the machinery developed in this work. More
interestingly, we will compute all the moduli spaces of deformations of
these submanifolds in codimension two where they naturally live, finding
the first known examples of honestly deformable Euclidean submanifolds
of rank and codimension two.


\vskip 0.3cm

Nothing similar to the local results for hypersurfaces due to Sbrana and
Cartan had been known for codimension higher than one until recently.
Locally, rank two elliptic and parabolic submanifolds were shown to be
honestly rigid in \cite{df3}, and it was not clear what would happen for
hyperbolic submanifolds. We will show here that, in contrast, hyperbolic
submanifolds are not genuinely rigid.

Now, for global rigidity, \cite{dm2} is devoted to show that any pair of
isometric immersions in codimension two of a compact Riemannian manifold
$M^n$ is nowhere genuine, provided certain mild singular extensions are
allowed. Indeed, the authors had to consider hyperbolic Sbrana-Cartan
hypersurfaces $N^{n+1}\subset\R^{n+2}$ together with their singular set,
$\Sigma^n\subset N^{n+1}$, which is itself, in fact, a regular deformable
hyperbolic submanifold in $\R^{n+2}$.
They called them {\it generalized Sbrana-Cartan hypersurfaces}. Actually,
when deforming $N^{n+1}$ its singular set $\Sigma^n$ also deforms, and
$M^n$ could very well share an open subset with $\Sigma^n$. More in the
spirit of this work, we can think of $\Sigma^n$ as a regular submanifold
in $\R^{n+2}$ that extends as a Sbrana-Cartan hypersurface
$N^{n+1}\subset \R^{n+2}$, but allowing singularities along $\Sigma^n$.
Yet, the actual necessity of considering singular extensions was not
established in \cite{dm2}, because the pair of immersions could also
extend regularly. Precisely this kind of singularities will appear again
in our \tref{genus} below, and as a consequence we will conclude that it
is indeed necessary to consider singular extensions to obtain this
strong genuine rigidity in codimension two for compact manifolds;
cf. \cref{neces}.

\begin{remark}\po\label{false} {\rm
We point out that the main result in \cite{dm2} is not actually correct
as it is stated in that paper. Indeed, by the same reason discussed
above, singular flat extensions may appear when the submanifold has
an open subset of flat points with first normal spaces of dimension less
than two. The omission is due to some minor gaps that have been
discovered recently by Felippe Guimar\~aes in some lemmas.
However, all gaps are fixed once singular flat extensions
are allowed, hence the result remains valid if we consider the flat
extensions together with their singularities. In any case, these gaps do
not affect our work since we will use these lemmas only when applied to
nonflat submanifolds with two dimensional first normal spaces.}
\end{remark}

\subsection{Shared dimension of a pair of curves}  

For later use, we introduce an elementary property
about curves in Euclidean space. Given two curves $\alpha_1$ and
$\alpha_2$ in $\R^N$ we define the
{\it shared dimension between $\alpha_1$ and $\alpha_2$},
denoted by $\overline I(\alpha_1,\alpha_2)$, as the smallest integer $k$
for which there is an orthogonal decomposition in affine subspaces,
$\R^N=\V_1\oplus \V^k \oplus \V_2$, satisfying
$\spa(\alpha_i)\subset \V_i\oplus \V^k$, $i=1,2$, where $\spa(\alpha)$
stands for the smallest affine linear subspace which contains the image
of $\alpha$. We call $\V^k$ the {\it shared subspace between $\alpha_1$
and $\alpha_2$}. Of course, generically,
$\overline I(\alpha_1,\alpha_2)=N$, and
$\overline I(\alpha_1,\alpha_2)=0$ if and only if the surface
$g(u,v)=\alpha_1(u)+\alpha_2(v)$ has flat normal bundle,
when $\alpha_1',\alpha_1'',\alpha_2',\alpha_2''$ are independent.
We need the following elementary characterization of
$\overline I(\alpha_1,\alpha_2)$.

\begin{lemma}\po\label{cos}
Given two curves $\alpha_1(u)$ and $\alpha_2(v)$ in $\R^N$, the integer
$\overline I(\alpha_1,\alpha_2)$ agrees with the minimum integer $k$ such
that $\la\alpha_1'(u),\alpha_2'(v)\ra$ can be written as a sum of the
form $\sum_{i=1}^ka_i(u)b_i(v)$ for certain smooth functions of one
variable $a_i,b_i$, $1\leq i\leq k$.
\end{lemma}
\proof
It is clear by definition that $k\leq \overline I(\alpha_1,\alpha_2)$. To
prove the opposite inequality, define
$\tilde\alpha_1=(\alpha_1, \int_0^u a_1(s)ds,\dots, \int_0^u a_k(s)ds)$
and
$\tilde\alpha_2=(\alpha_2,-\int_0^v b_1(s)ds,\dots,-\int_0^v b_k(s)ds)$
as orthogonal curves in $\R^{N+k}$. So,
$\R^{N+k}=\tilde\V_1^{n_1}\oplus^\perp \tilde \V_2^{n_2}$ with
$\spa(\tilde\alpha_i)\subset\tilde \V_i^{n_i}$, $i=1,2$.
Consider $\V_i=\tilde \V_i^{n_i}\cap(\R^N\!\times\!\{0\})\subset\R^N$,
and complete to an orthogonal decomposition,
$\R^N=\V_1\oplus\V\oplus\V_2$. By construction, $\spa(\alpha_i)$ is
orthogonal to $\V_j$, $1\leq j\neq i\leq 2$, and
$\dim \V = N-\dim \V_1-\dim \V_2 \leq N-(n_1-k)-(n_2-k)=k$.
\qed

\medskip

We will actually need the local version of this concept. Since
$\overline I$ does not increase when we restrict the domains of the
curves, we define the {\it local shared dimension between $\alpha_1$ and
$\alpha_2$} as the integer-valued function
$$
I(\alpha_1,\alpha_2)(u,v):=\lim_{\epsilon\to 0}
\overline I\left(\alpha_1|_{(u-\epsilon,u+\epsilon)},
\alpha_2|_{(v-\epsilon,v+\epsilon)}\right).
$$
Since $I(\alpha_1,\alpha_2)$ is clearly semicontinuous, it
is constant along connected components of an open dense subset
of the parameters $(u,v)$.

\section{Projecting the data}\label{project} 

In this section we show that all the data relevant to the study of the
rigidity of a hyperbolic submanifold actually project to the leaf space
of the relative nullity foliation, $L^2=M^n/\Delta$. Notice that $L^2$ is
a smooth surface when working locally on $M^n$, while globally it may
only fail to be Hausdorff. The objects and notations introduced in this
section will be used throughout the whole work.

\medskip

Let $f:M^n\rightarrow\R^{n+2}$ be a hyperbolic submanifold.
We denote by $\{ \xi_1, \xi_2 \} \subset T^\perp_fM$ and
$\{ Y_1,Y_2 \} \subset \Delta^\perp$ local smooth unit frames such that
\begin{equation}\label{ker1g}
\ker A_{\xi_i}=\Delta\oplus\spa\{Y_i\},\ \ \ i=1,2,
\end{equation}
where we also require that the nonzero eigenvalue $\lambda_i$ of
$A_{\xi_i}$ is positive. Let $0<\theta<\pi$ and $0<\omega<\pi$ be the
angles between $\xi_1$ and $\xi_2$, and between $Y_1$ and $Y_2$,
respectively. The function $\theta$, which we call the {\it main angle}
of $f$, plays a crucial role in this work. Fix the orientation on
$T^\perp_fM$ and $\Delta^\perp$ determined by these bases, and complete
to local smooth oriented orthonormal frames
$\{\xi_i,\eta_i\}\subset T^\perp_fM$ and
$\{X_i,Y_i\} \subset\Delta^\perp$, i.e.,
\begin{equation}\label{chbag}
\sin(w)X_i = (-1)^{i+1} (\cos(w)Y_i-Y_j), \ \
\sin(\theta)\eta_i = (-1)^i (\cos(\theta)\xi_i-\xi_j), 
\end{equation}
$1\leq i\neq j\leq 2$.
Observe that, since we are in codimension two and rk $A_{\xi_1}=1$,
Gauss equation reduces to ${\rm scal}=\det (A_{\eta_i}|_{\Delta^\perp})$, 
where scal denotes the (non-normalized) scalar curvature of $M^n$.
Then from \eqref{ker1g} and \eqref{chbag} we obtain
\begin{equation}\label{gaussg}
{\rm scal}=
-\sin^2(\omega)\,\frac{\cos(\theta)}{\sin^2(\theta)}\lambda_1\lambda_2.
\end{equation}

Finally, denote by $\psi^i$ the normal connection 1-form associated to
the frame $\{\xi_i,\eta_i\}$,
$$
\psi^i(X)=\la\nabla^\perp_X\xi_i,\eta_i\ra.
$$
Although the normal connection 1-form depends on the chosen
orthonormal frame, its differential is independent
since by \eqref{chbag} we have
\begin{equation}\label{nc1fg}
\psi^2=\psi^1+d\theta.
\end{equation}

\medskip

We proceed to obtain information about the evolution of the geometric
data described above along the relative nullity foliation $\Delta$.
As a consequence, we will be able to project
all the data relevant to the understanding of the isometric
deformation problem to its leaf space
$$
\pi:M^n\rightarrow L^2:=M^n/\Delta.
$$
Recall that $L^2$ is intrinsic when
$M^n$ is nowhere flat, since it coincides with $M^n/\Gamma$.
\begin{proposition}\po\label{tptqg}
With the above notations, the following holds for $1\leq i\neq j\leq 2$:
\begin{enumerate}[$i)$]
\item $\xi_i$ is constant in $\R^{n+2}$ along the leaves of $\Delta$;
\item $Y_i$ is constant in $\R^{n+2}$ along the leaves of $\Delta$;
\item $\{Y_1,Y_2\}$ is an eigenbasis of the hyperbolic structure
on $\Delta^\perp$, and the splitting tensor of $\Delta$ satisfies that
$C_TY_i=T(\ln\lambda_j)Y_i,\ $ for all $\ T\in\Delta$;
\item There exist functions $\overline{\theta}$,
$\overline{\omega}$ and a 1-form $\overline{\psi^i}$
on $L^2$ such that $\theta=\overline{\theta}\circ\pi$,
$\omega=\overline{\omega}\circ\pi$ and
$\psi^i=\pi^*\overline{\psi^i}$.
\end{enumerate}
\end{proposition}
\proof
The Codazzi equation for $(A_{\xi_i},T\in\Delta,Z\in\Delta^\perp)$ yields
\begin{equation} \label{codag}
\nabla_TA_{\xi_i}Z-A_{\xi_i}[T,Z]-\psi^i(T)A_{\eta_i}Z=0.
\end{equation}
This for $Z=Y_i\in\ker A_{\xi_i}$ gives
\begin{eqnarray}\label{4}
A_{\xi_i}[T,Y_i]=-\psi^i(T)A_{\eta_i}Y_i.
\end{eqnarray}
Since $Y_i\in\ker A_{\xi_i}$ and $f$ has rank two, we have that
$0\neq A_{\eta_i} Y_i\in\im A_{\xi_j}$.
However, since $\im A_{\xi_1}\cap\im A_{\xi_2}=0$, both sides of
\eqref{4} vanish, and then
\begin{equation}\label{coda9g}
\psi^i(T)=0,\ \ \ \forall\, T\in\Delta,
\end{equation}
\vskip -0.5cm
\begin{equation}\label{coda1g}
[T,Y_i]\in\ker A_{\xi_i}\ \ \ \forall\, T\in\Delta.
\end{equation}
Equation \eqref{coda9g} says that $\xi_i$ is parallel along $\Delta$
with respect to the normal connection, which by the Weingarten formula
implies $(i)$.

By the above, \eqref{codag} reduces to
\begin{equation}\label{coda2g}
\nabla_TA_{\xi_i}Z=A_{\xi_i}[T,Z],
\end{equation}
which shows that the line bundle $\im A_{\xi_i}$ is parallel along
$\Delta$ with respect to the Levi-Civita connection $\nabla$ of $M^n$. In
particular, the unitary vector field $X_i$, and then also $Y_i$, are 
parallel along $\Delta$ with respect to $\nabla$, which proves $(ii)$ by
the Gauss formula.

Now, $(ii)$ says that $[T,Y_i]=-\nabla_{Y_i}T$. Taking the orthogonal
projection onto $\Delta^\perp$ of this relation we conclude from
\eqref{coda1g} that $\{Y_1,Y_2\}$ is a common eigenbasis for all
splitting tensors and thus, by \eqref{stdhsg}, for the hyperbolic
structure as well since $f$ is nowhere surface-like. In other words,
there are 1-forms $b_j$ on $\Delta$ such that
\begin{equation}\label{evstg}
C_TY_j=b_j(T)Y_j, \ \ \ \forall\, T\in\Delta.
\end{equation}
A straightforward computation using \eqref{chbag} gives
$$
b_j(T)=\langle C_TX_i,X_i\rangle.
$$
On the other hand, setting $Z=X_i$ in \eqref{coda2g} and using
$A_{\xi_i}X_i=\lambda_iX_i$ we obtain that the right-hand side of the above
equation equals $T(\ln\lambda_i)$, from which $(iii)$ follows.

By $(i)$ and $(ii)$ and the definition of $\theta$ and $\omega$
both angles are constant along $\Delta$, so they project to $L^2$.
Finally, the Ricci equation for
$(\xi_i,\eta_i,T\in\Delta,Z\in\Delta^\perp)$ yields
$$
d\psi^i(T,Z)=0,\ \forall T\in\Delta,\ Z\in\Delta^\perp,
$$
which alongside \eqref{coda9g} implies that $\psi^i$ is projectable by
Corollary 12 in \cite{dft3}.
\qed
\medskip

We often identify projectable functions and 1-forms with their
respective projections, when there is no risk of confusion.
This will be further clarified in the next section after \tref{param}.

\medskip

We show next that $(ii)$ and $(iii)$ above allow us to rescale $Y_1$ and
$Y_2$ so that the resulting frame projects to a coordinate frame on
$L^2$. These coordinates will be used throughout this work, and in
particular will be useful to prove the existence of polar surfaces.

\begin{proposition}\po\label{cspg}
There exist smooth positive functions $\mu_1$ and $\mu_2$ on $M^n$ and a
coordinate system $(u_1,u_2)$ on $L^2$ such that the frame $\{Z_1,Z_2\}$
defined by $Z_i=\mu_iY_i$ satisfies
$$\partial_{u_i}\circ\pi=\pi_*\circ Z_i,\ \ i=1,2.$$
\end{proposition}
\proof
According to Proposition 10 in \cite{dft3}, the necessary and sufficient
condition for the vector fields $Z_i$ to be projectable is that
\begin{equation}\label{prong}
[Z_i,T]\in\Delta,\ \ \forall T\in\Delta,
\end{equation}
whereas the projections $\pi_*\circ Z_i$ come from a local coordinate
system if, additionally,
\begin{equation}\label{procg}
[Z_1,Z_2]\in\Delta.
\end{equation}
 From \pref{tptqg}-$(ii)$ and $(iii)$, we have that \eqref{prong} reduces
to
$T(\mu_i)=-b_i(T)\mu_i$,
with $b_i=d(\ln\lambda_j)|_\Delta$, $1\leq i\neq j\leq2$,
while condition \eqref{procg} can be written as
$Y_j(\mu_i)=-r_i \mu_i$,
for $r_1$ and $r_2$ defined by $[Y_1,Y_2]+r_1Y_1-r_2Y_2\in\Delta$.
In other words, we must show that the first order system of PDEs
\begin{equation}\label{pdesg}
d(\ln\mu_i)|_\Delta = -b_i\,,\ \ \
Y_j(\ln\mu_i) = -r_i
\end{equation}
is integrable.

Consider the distribution $\Omega_i=\Delta\oplus\spa\{Y_j\}$ on $M^n$,
$1\leq i\neq j\leq2$. Since $Y_j$ is parallel along $\Delta$ and an
eigenvector of all $C_T$, $T\in\Delta$, we have from the integrability
of~$\Delta$ that $\Omega_i$ is also integrable. Define on $\Omega_i$ the
1-form $\sigma_i$ by
$$
\sigma_i|_\Delta=-b_i,\ \ \sigma_i(Y_j)=-r_i.
$$
Since all our considerations are local, the integrability condition of
\eqref{pdesg} translates into the exactness of $\sigma_i$.
Since the 1-form $b_i$ on $\Delta$ is exact,
$d\sigma_i|_{\Delta\times\Delta}=0$. Thus, it suffices to show that
$d\sigma_i(T,Y_j)=0$, or equivalently,
\begin{equation}\label{intcg}
Y_j(b_i(T))-T(r_i)=b_i(\nabla^v_{Y_j}T)-r_ib_j(T),
\ \ \ \forall\, T\in\Delta,
\end{equation}
since, by \pref{tptqg}-$(ii)$,
$$\sigma_i([T,Y_j])=-\sigma_i(\nabla_{Y_j}T)
=\sigma_i(C_TY_j)-\sigma_i(\nabla^v_{Y_j}T)
=b_i(\nabla^v_{Y_j}T)-r_ib_j(T),$$
where $\nabla^h$ and $\nabla^v$ stand for the connections induced by
$\nabla$ on $\Delta^\perp$ and $\Delta$, respectively.

It is well-known and easy to see that the Codazzi equation for $f$
implies that the splitting tensor $C$ itself is a Codazzi tensor, that is,
$(\nabla^h_{Y_i}C_T){Y_j}-C_{\nabla^v_{Y_i}T}{Y_j}=
(\nabla^h_{Y_j}C_T){Y_i}-C_{\nabla^v_{Y_j}T}{Y_i}$.
This can be easily rewritten using \eqref{evstg} and \eqref{chbag} as
\begin{equation}\label{niu}
Y_j(b_i(T))-b_i(\nabla^v_{Y_j}T)=
(-1)^i\sin(\omega)^{-1}(b_i(T)-b_j(T))G_i,
\end{equation}
where $G_i=\langle \nabla_{Y_i}Y_j,X_j\rangle$.
By \eqref{chbag} we express $r_i$ in terms of $G_1$ and $G_2$ as
\begin{equation}\label{rrr}
r_i=(-1)^i\sin(\omega)^{-1}(G_i-\cos(\omega)G_j).
\end{equation}
Now, as $T\in\Delta=\Gamma$, we have from \pref{tptqg}-$(ii)$ that
\begin{equation}\label{ctn1g}
0=\langle R(T,Y_i)Y_j,X_j\rangle
=\langle\nabla_T\nabla_{Y_i}Y_j,X_j\rangle-b_i(T)G_i=T(G_i)-b_i(T)G_i.
\end{equation}
Since $\omega$ and $X_i$ are constant along $\Delta$, differentiate
\eqref{rrr} along $\Delta$ and use \eqref{ctn1g} to conclude that
$T(r_i)-r_ib_j(T)$ agrees with the right-hand side of \eqref{niu}.
This proves \eqref{intcg}, as wished.

Therefore, the system \eqref{pdesg} is integrable, which means that
$\mu_i$ can be arbitrarily prescribed along a fixed integral curve
$\gamma$ of $Y_i$ and then extended along each leaf of $\Omega_i$ through
$\gamma$ as a solution of \eqref{pdesg}.
\qed

\section{Polar surfaces and the parametrization} \label{psatpg}  
In this section we show that the normal space of a hyperbolic submanifold
is always integrable via a so-called polar surface. We use this to
recover any hyperbolic submanifold from its polar surface through a very
simple and explicit parametrization. As an important application, we will
classify all flat Euclidean submanifolds in codimension two in a simple
and explicit way.

\medskip

Given a hyperbolic submanifold $f:M^n\to\R^{n+2}$, we will show
first the (local) existence of a surface $g:L^2=M^n/\Delta\to\R^{n+2}$
that integrates its normal space in the sense that
$T^\perp_fM=\pi^*(T_gL)$. We follow all the notations and definitions of
the previous section and, for convenience and from now on, $u=u_1$ and
$v=u_2$ as subindexes will denote the corresponding partial derivatives
for the local coordinates $(u,v)$ constructed in \pref{cspg}.

Assume that $L^2$ is simply-connected, and suppose that there is such a
surface $g$. Then, $g_u=c\,\overline\xi_1+a\,\overline\xi_2$ and
$g_v=b\,\overline\xi_1+d\,\overline\xi_2$, for certain functions
$a,b,c,d$ over $L^2$. Differentiating the first equation with respect to
$v$ and the second with respect to $u$, and projecting orthogonally onto
$T^\perp_{g(\pi(x))}L\supseteq\Delta^\perp_f(x)$ we get
$cA_{\xi_1}Z_2=dA_{\xi_2}Z_1$. Then $c=d=0$ since $A_{\xi_1}Z_2$ and
$A_{\xi_2}Z_1$ are linearly independent.
So consider the following first order system of PDEs over $L^2$:
\begin{equation}\label{sysg}
g_u = a\ \overline\xi_2\,,\ \ \ \
g_v = b\ \overline\xi_1.
\end{equation}
Differentiating the first equation in \eqref{sysg} with respect to $v$
gives
$$
g_{uv}\circ \pi = a_v\circ\pi\ \xi_2+a\circ\pi\ 
\tilde \nabla_{\pi_*\circ Z_2}\overline\xi_2 =
a_v\circ \pi\ \xi_2 + a\circ \pi\ \tilde\nabla_{Z_2}\xi_2.
$$
Since $Z_2\in\ker A_{\xi_2}$ and the normal connection form
$\psi^2=\la\nabla^\perp_{\, \bullet}\xi_2,\eta_2\ra$ projects to
$\overline{\psi^2}$, we get
\begin{equation}\label{2v}
g_{uv}=
a_v\,\overline\xi_2+a\,\overline{\psi^2}(\partial_v)\overline\eta_2,
\end{equation}
and analogously for the second equation in \eqref{sysg},
\begin{equation}\label{1u}
g_{vu}=
b_u\,\overline\xi_1+b\,\overline{\psi^1}(\partial_u)\overline\eta_1.
\end{equation}
Since both bases are equally oriented, we get that the integrability
conditions for \eqref{sysg} are
$$
\sin(\overline\theta)\,a_v=\overline{\psi^1}(\partial_u)\,b
-\cos(\overline\theta)\overline{\psi^2}(\partial_v)\,a,
\ \ \
\sin(\overline\theta)\,b_u=\cos(\overline\theta)
\overline{\psi^1}(\partial_u)\,b-\overline{\psi^2}(\partial_v)\,a.
$$
 From $\sin(\overline\theta)>0$ we conclude that this system always has
solutions $a\neq 0,b\neq 0$ and so by \eqref{sysg} there exists a regular
surface $g:L^2\to\R^{n+2}$, with
\begin{equation}\label{polar}
\pi^*(T_{g}L) = T^\perp_fM,\ \ \ 
{\rm and}\ \ \ \pi^*(N^1_g)=\Delta^\perp_f.
\end{equation}
In addition,
\begin{equation}\label{fff}
E=a^2,\ \ G=b^2,\ \ F=ab\,\cos(\overline\theta)
\end{equation}
are the coefficients of the first fundamental form of $g$. Moreover, $g$
is also hyperbolic in the sense that $\dim N^1_g=2$ and
$\alpha_g(\partial_u,\partial_v)=0$ for the second fundamental form
$\alpha_g$ of $g$, since $g_{uv}\in T_gL$. Recall that a
coordinate system $(u,v)$ on $L^2$ such that
$\alpha_g(\partial_u,\partial_v)=0$ everywhere is called {\it conjugate}.
Following \cite{df1} we call such a surface $g$ a {\it polar surface
of~$f$}, and from now on we consider on $L^2$ the metric induced by a
(fixed) polar surface $g$ of $f$, as in \eqref{fff}.

\begin{remark}\po\label{curvrels}
{\rm The 1-forms $\overline{\psi^i}$ are tangent connection forms for
$g$, so $d \overline{\psi^i} = K_g\ dA$,
where $K_g$ denotes the Gaussian curvature of $g$ and $dA$ its area
element. Since $\psi^i=\pi^*\overline{\psi^i}$, we have for the normal
curvature 2-form $R^\perp_f=d\psi^i$ of $f$ that
$$
R^\perp_f = (K_g\circ\pi)\ \pi^*dA.
$$
In particular, $g$ is flat if and only if $f$ has flat normal bundle.
Furthermore, since $\theta$ is the angle between the conjugate
directions of $g$, we also have that $f$ is flat if and only if $g$ has
flat normal bundle, that is, $\theta\equiv\pi/2$, or, equivalently,
$F\equiv 0$.}
\end{remark}

\begin{remark}\po\label{trivial}
{\rm Hyperbolic surfaces in $\R^k$ are trivial to construct and classify
locally: they are simply given by $k$ generic solutions of a fixed
wave equation. Indeed, let $U\subset\R^2$ be an open set with coordinates
$(u,v)$, and let $\Gamma^u,\Gamma^v\!:U\to\R$ be two arbitrary smooth
functions. Consider the second order linear wave differential operator
\begin{equation}\label{q}
Q = \partial_u\circ\partial_v-\Gamma^u\partial_u-\Gamma^v\partial_v.
\end{equation}
Take $k$ smooth functions $g=(g_1,\dots,g_k)$ that are solutions of the
wave equation $Q=0$, which are {\it independent} in the sense that
$g_u,g_v,g_{uu},g_{vv}$ are pointwise linearly independent. Then,
$g:U\to\R^k$ is a regular surface, $\dim N^1_g=2$, and
\begin{equation}\label{guv}
g_{uv}-\Gamma^ug_u-\Gamma^vg_v=0,
\end{equation}
thus $g$ is hyperbolic. In this sense, there is no geometry involved in
the construction of hyperbolic surfaces. Notice also that
$\Gamma^u,\Gamma^v$ are automatically Christoffel symbols of $g$ since
$\nabla^g_{\partial_u}\partial_v=\Gamma^u\partial_u+\Gamma^v\partial_v$.
In particular, $Q$ is related with the Hessian of the metric induced
by $g$ by the relation $Q(\rho)=\Hess_\rho(\partial_u,\partial_v)$.}
\end{remark}

\begin{remark}\po\label{wd}
{\rm By \eqref{chbag}, \eqref{2v}, \eqref{1u}, \eqref{fff}
and \eqref{guv} it holds that
\begin{equation}\label{2v1u}
\psi^1(\partial_u)=\sin(\theta)\sqrt{E/G}\ \Gamma^u,\ \ \ {\rm and}\ \ \
\psi^2(\partial_v)=-\sin(\theta)\sqrt{G/E}\ \Gamma^v.
\end{equation}
Using \eqref{nc1fg} we see that these two equations recover the normal
connection of $f$ in terms of the metric of $g$. In particular, although
a polar surface of a given hyperbolic submanifold is not unique, the
functions $F/\sqrt{EG}=\cos(\theta)$ as well as $\sqrt{E/G}\ \Gamma^u$
and $\sqrt{G/E}\ \Gamma^v$ coincide for all its polar surfaces.}
\end{remark}


Next we proceed to show how to recover any hyperbolic submanifold $f$
from a polar surface $g$ of it by providing a simple parametrization
of $f$ that depends only on $g$ and a smooth function over $L^2$
satisfying the same wave equation as $g$.

\medskip

Consider $h:L^2=M^n/\Delta\to M^n$ any local {\it cross-section} to the
relative nullity foliation, that is, $h^*(\Delta)\oplus T_hL = h^*(TM)$.
Since the leaves of relative nullity are mapped by $f$ to (open subsets
of) affine $(n-2)$-dimensional subspaces, \eqref{polar} implies that
$$
\Psi:(N^1_g)^\perp\subset T_g^\perp L\to\R^{n+2},\ \ \ \ 
\Psi(\mu) = (f\circ h) + \mu
$$
for any smooth section $\mu$ of $(N^1_g)^\perp$ parametrizes, at regular
points, our submanifold $f$. Notice that, with the notations of
\rref{trivial},
$$
(N^1_g)^\perp = \spa\{g_u,g_v,g_{uu},g_{vv}\}^\perp.
$$
We will conclude the parametrization by describing $h':=f\circ h$ in
terms of $g$ and a smooth function over $L^2$.

By \eqref{polar}, the map $h'$ is characterized by the property
\begin{equation}\label{cs}
f_*(h^*(\Delta))\oplus T_{h'}L=h^*(T_fM)=T^\perp_gL.
\end{equation}
Decompose the position vector $h'\in \R^{n+2}=T_gL\oplus T_g^\perp
L$. The component in $(N^1_g)^\perp$ of $h'$ plays no role in
\eqref{cs} (or, equivalently, it can be absorbed in the parametrization
$\Psi$), so we write
$$
h'=X+\eta,\ \ X\in T_gL,\ \ \eta\in N^1_g.
$$
But with this decomposition \eqref{cs} is equivalent to $\nabla^g_\bullet
X = A_\eta^g$. In particular, $X$ is the gradient of a certain function
$\rho$, $X=\nabla\rho$, whose Hessian satisfies
\begin{equation}\label{hess}
\Hess_\rho = A_\eta^g.
\end{equation}
Since $\dim N^1_g =2$, by dimension reasons, given a function $\rho$
on $L$ there exists (a unique) $\eta\in N^1_g$ for which \eqref{hess}
holds if and only if $\Hess_\rho(\partial_u,\partial_v)=0$. In other
words, by \rref{trivial}, $\rho$ satisfies the same linear PDE as the
coordinate functions of $g$, i.e., \mbox{$Q(\rho)=0$} for $Q$ as in
\eqref{q}.

\medskip

This was the final ingredient for the main result in this section,
that is interesting in its own right:

\begin{theorem}\po\label{param}
Let $g=(g_1,\dots,g_{n+2}):W\subset\R^2\to\R^{n+2}$ be $(n+2)$
independent solutions of a wave equation \eqref{q}, $Q(g)=0$, and let
$\rho$ be another solution, $Q(\rho)=0$. Then, $g$ is a hyperbolic
surface once we endow $W$ with the metric induced by $g$, and there is a
unique $\eta_\rho\in N^1_g$ such that $\Hess_\rho=A^g_{\eta_\rho}$.
Moreover, the map
\begin{equation}\label{paramet}
\Psi:(N^1_g)^\perp\subset T_g^\perp L\to\R^{n+2},\ \
\Psi(\mu_x) = (\nabla\rho + \eta_\rho)(x) + \mu_x,
\end{equation}
parametrizes, at regular points, a hyperbolic $n$-dimensional
submanifold of $\R^{n+2}$ for which $g$ is a polar surface.

Conversely, any hyperbolic Euclidean submanifold in codimension two
can be locally parametrized in this way.
\end{theorem}
\proof
We have already proved the converse claim. For the direct statement,
first notice that, by construction, the submanifold $\Psi$ has $g$ as a
polar surface at its regular points, i.e., $\pi^*(T_gL) = T^\perp_\Psi M$
for $\pi:M^n=(N^1_g)^\perp\to L^2$ and $M^n$ endowed with the metric
induced by $\Psi$. In fact, since $g$ is constructed from independent
solutions, $\dim N^1_\Psi=2$. In addition, $g_u$ and $g_v$ are
independent normal vector fields to $\Psi$ whose shape operators have
rank one since
$(g_v)_u \circ \pi = (g_u)_v \circ \pi
= g_{uv} \circ \pi \in \pi^*(T_gL) = T^\perp_\Psi M$.
Therefore, $\Psi$ is hyperbolic.
\qed
\vspace{1.5ex}

In a local trivialization $(u,v,t_1,\dots,t_{n-2})$ of
$(N^1_g)^\perp$ determined by a moving frame
$\{\eta_1,\,\dots,\,\eta_{n-2}\}$ of $(N^1_g)^\perp$ along a
conjugate local coordinate system $(u,v)$ of $L^2$, $\Psi$ can be
written as
$$
\Psi(u,v,t_1,\dots,t_{n-2})
=(\nabla\rho + \eta_\rho)(u,v)+\sum_{i=1}^{n-2}t_i\eta_i(u,v).
$$
We will consider from now on this as the {\it standard coordinate system}
of our hyperbolic submanifolds. In particular, the coordinate vector
fields $\partial_u$, $\partial_v$ as well as differentiation with respect
to $u$ and $v$ makes sense now also on $M^n$, and the fact that $u$ and
$v$ are considered as coordinate functions in both $M^n$ and $L^2$ will
not cause any confusion. For example, since $(u,v)$ are conjugate
coordinates for~$g$, the coordinate vector fields $\partial_u$ and
$\partial_v$ are also conjugate for $\Psi$, i.e.,
$\alpha_\Psi(\partial_u,\partial_v)=0$. Moreover, that a certain function
$h=h(u,v,t_1,\dots,t_{n-2})$ on $M^n$ projects to its leaf space
$L^2$ means simply that it does not depend on the coordinates
$t_1,\dots,t_{n-2}$. So we will always denote with the same symbol
$h=h(u,v)$ the projection of $h$ to $L^2$.

\begin{remark}\po
{\rm
 From \tref{param} it follows that $f$ splits a Euclidean factor if and
only if $g$ is not substantial, that is, if $g$ reduces codimension.
Moreover, it is easy to check that $f$ is surface-like if and only if $g$
has substantial conformal codimension 2, i.e., the image of $g$ is
contained in a 4-dimensional umbilical submanifold of the ambient space.}
\end{remark}

As an immediate consequence of \tref{param} and \rref{curvrels} we have a
parametric description of all generic rank two Euclidean submanifolds $f$
in codimension two with flat normal bundle in terms of Euclidean
hyperbolic flat surfaces. Generic here means simply that \mbox{$\dim
N^1_f=2$} everywhere, a condition that, by \pref{n12}, is automatic
if~$N^1_f$ has constant dimension and $f$ is not a composition. More
interestingly, in a dual way and also by \rref{curvrels}, we can
characterize all generic flat Euclidean submanifolds in codimension two:

\begin{corollary}\po\label{flat}
Let $g:L^2\to\R^{n+2}$ be any surface with flat normal bundle and
principal coordinates $(u,v)$, and $\rho:L^2\to\R$ any smooth function
satisfying $\Hess_\rho(\partial_u,\partial_v)=0$. Then, the map
\eqref{paramet} parametrizes, at regular points, a flat $n$-dimensional
submanifold of $\R^{n+2}$ for which $g$ is a polar surface.
Conversely, any generic flat submanifold in codimension two
can be locally parametrized this way.
\end{corollary}

The construction in \cref{flat} is much simpler, direct and explicit than
the one given in Theorem 13 in \cite{df0}. While the latter depends on an
elusive kind of surfaces in the sphere, called of `type $C$', any surface
in $\R^k$ with flat normal bundle parametrized by lines of curvature can
be described explicitly from $k$ arbitrary solutions of a certain simple
linear integrable system of PDEs thanks to the beautiful construction due
to E. Ferapontov in \cite{fera} (see also \cite{dft1} for the
generalization of this construction to arbitrary dimensions). In
particular, \cref{flat} tells us what all those surfaces $\xi$ of
type~$C$ really are: they are just the principal directions of any
surface $g$ with flat normal bundle, i.e., $\xi =
g_u/\|g_u\|:L^2\to\Sp^{n+1}$.

\begin{remark}\po\label{comp2}
{\rm As already pointed out, one of the main reasons to discard
compositions when studying Euclidean rigidity in codimension two is the
existence of a very simple local classification of flat Euclidean
hypersurfaces (via the Gauss parametrization). Now \cref{flat} also
justifies discarding compositions when working in codimension~3, as well
as the introduction of our concept of honest rigidity.}
\end{remark}

\section{The 6 invariants of a hyperbolic submanifold}\label{ihs}  
In this section we extract the main data of a Euclidean
hyperbolic submanifold in codimension two and show how it completely
determines its deformations by means of a set of six
invariants.

\medskip

Let $f:M^n\to\R^{n+2}$ be a nowhere flat hyperbolic
submanifold, and $L^2=M^n/\Delta$ its nullity leaf space.
Following the notations of \sref{project}, we set
$s:=\sin^2(\theta)$ for the main angle $\theta$ of $f$, and, in a
standard coordinate system, we define the {\it main symbols} of $f$ by
\begin{equation}\label{lambdasdefs}
\Lambda^u:=\frac{\psi^2(\partial_u)}{\tan(\theta)},\ \ \
\Lambda^v:=-\frac{\psi^1(\partial_v)}{\tan(\theta)}.
\end{equation}
We can easily express $\Lambda^u$, $\Lambda^v$ and $s$, which by
\pref{tptqg}-$(iv)$ can and will be seen as functions over $L^2$, in
terms of the first fundamental form of the polar surface $g$. Indeed,
from \eqref{nc1fg} and \eqref{2v1u} we get
\begin{equation}\label{lambdas}
\Lambda^u=\frac{F}{G}\,\Gamma^u+\frac{s_u}{2s},\ \ \
\Lambda^v=\frac{F}{E}\,\Gamma^v+\frac{s_v}{2s},\ \ \
s=1-\frac{F^2}{EG},
\end{equation}
where $E$,$F$ and $G$ are the coefficients of the first fundamental form
of $g$ in the coordinate system $(u,v)$. Notice that, by \rref{wd}, the
expressions in \eqref{lambdas} do not depend on the particular choice of
a polar surface of $f$. We set
$$
\D(f):=\{\theta,\Lambda^u,\Lambda^v,\kappa^u,\kappa^v\},
$$
where
$$
\kappa^u:=\lambda_1^2/s,\ \ \ \kappa^v:=\lambda_2^2/s.
$$
Observe that by \eqref{nc1fg} and
$\lambda_1,\lambda_2>0$, we can recover $\lambda_1,\,\lambda_2,\,\psi^1$
and thus the second fundamental form and normal connection of $f$ from
$\D(f)$. The Fundamental Theorem of Submanifolds hence says that
$\D(f)$ determines $f$ itself when $M^n$ is simply-connected.

\medskip

Our following result provides a set of six functions that both determine
and remain invariant under deformations $\hat f$ of $f$ that are nowhere
compositions. Observe that this is a slightly different approach to the
usual ones that deal with the study of deformations. From now on, we add
a hat to indicate the objects of $\hat f$ corresponding to those of $f$.
In particular, $\hat \theta$ refers to the main angle of $\hat f$, and
$\hat \Lambda^u,\hat \Lambda^v$ to its main symbols.

\begin{proposition}\po\label{ipg}
Let $f:M^n\rightarrow\R^{n+2}$ be a nowhere flat hyperbolic submanifold.
Then, the functions ${\cal G}$,
$\CC^u_1$, $\CC^v_1$, $\CC^u_2$, $\CC^v_2$ and $\RR$ given by
\begin{equation}\label{isg}
\begin{gathered}
{\cal G}=\cos(\theta)\sqrt{\kappa^u\kappa^v},\\
\CC^u_1=\kappa^u\Lambda^u,\ \ \ \CC^v_1=\kappa^v\Lambda^v,\\
\CC^u_2=\frac{\kappa^u_u}{\kappa^u}+2\Lambda^u,\ \ \ 
\CC^v_2=\frac{\kappa^v_v}{\kappa^v}+2\Lambda^v,\\
\RR=\rho_{uv}+\cos^2(\theta)\,\rho_u\rho_v-\Lambda^v\rho_u-\Lambda^u\rho_v,
\end{gathered}
\end{equation}
where $\rho=\ln(|\tan(\theta)|)$, are invariant, i.e., such
functions are preserved by any hyperbolic deformation of $f$.
Moreover, the ratios
$\tau^u:=\hat\kappa^u/\kappa^u$ and $\tau^v:=\hat\kappa^v/\kappa^v$
project to $L^2$.

Conversely, if $M^n$ is simply-connected, given functions
$\tau^u>0$, $\tau^v>0$, $0<\hat\theta<\pi$,
$\hat\Lambda^u$ and $\hat\Lambda^v$ on $L^2$ such that
$\hat\D=\{\hat{\theta}, \hat\Lambda^u,\hat\Lambda^v,
\hat\kappa^u:=\tau^u\kappa^u, \hat\kappa^v:=\tau^v\kappa^v\}$
satisfies system \eqref{isg}, there exists an isometric immersion
$\hat{f}:M^n\to\R^{n+2}$ with $\D(\hat{f})=\hat\D$.
\end{proposition}
\proof
Since $M^n$ is nowhere flat and $\hat f$ has rank two,
$\hat{\Delta}=\Gamma=\Delta$ and $\hat C=C$ everywhere. Moreover, since
$\hat{f}$ is nowhere surface-like, by \eqref{stdhsg} $\hat{f}$ is
hyperbolic with $\hat{J}=J$. In particular, by \pref{tptqg}-$(iii)$,
$\hat{Z}_i=Z_i$, $i=1,2$, and hence $\hat{\omega}=\omega$, so that the
invariance of ${\cal G}$ is simply a restatement of the Gauss equation
\eqref{gaussg} for $\Delta^\perp$ and $\hat f$. It also follows from
\pref{tptqg}-$(iii)$ that $\tau^u$ and $\tau^v$ project to $L^2$.

 From \eqref{ker1g} and \eqref{chbag}, the Codazzi
 equation for $(A_{\xi_1},\partial_u,\partial_v)$ is
\begin{equation}\label{codvg}
\sin(\theta)\na_{\partial_u}A_{\xi_1}\partial_v=
-\cos(\theta)\psi^1(\partial_u)A_{\xi_1}\partial_v
-\psi^1(\partial_v)A_{\xi_2}\partial_u.
\end{equation}
Setting $\delta_i=\hat{\lambda_i}/\lambda_i$, we have that
$\hat A_{\hat\xi_i}=\delta_iA_{\xi_i}$ and \eqref{codvg}
for $\hat{f}$ gives
\begin{equation}\label{codvhg}
\sin(\hat{\theta})((\delta_1)_uA_{\xi_1}\partial_v
+\delta_1\na_{\partial_u}A_{\xi_1}\partial_v)=
-\cos(\hat\theta)\hat\psi^1(\partial_u)\delta_1A_{\xi_1}\partial_v
-\hat\psi^1(\partial_v)\delta_2A_{\xi_2}\partial_u.
\end{equation}
Using \eqref{codvg} in \eqref{codvhg} and $\sin(\theta)\neq 0$ we get
$$
\left(\delta_1\sin(\hat\theta)\psi^1(\partial_v)-
\delta_2\sin(\theta)\hat\psi^1(\partial_v)\right)A_{\xi_2}\partial_u=
\hskip 4cm
$$
$$
\left((\delta_1)_u\sin(\theta)\sin(\hat\theta)+
\delta_1\sin(\theta)\cos(\hat\theta)\hat\psi^1(\partial_u)-
\delta_1\sin(\hat\theta)\cos(\theta)
\psi^1(\partial_u)\right)A_{\xi_1}\partial_v.
$$
Since $A_{\xi_1}\partial_v$ and $A_{\xi_2}\partial_u$ are linearly
independent everywhere, we get from the invariance of ${\cal G}$
the invariance of $\CC^v_1$ and 
$(\ln\delta_1)_u=\cot(\theta)\psi^1(\partial_u)
-\cot(\hat\theta)\hat\psi^1(\partial_u)$, which is clearly equivalent
to the invariance of $\CC^u_2$ by \eqref{nc1fg}. Similarly,
the invariance of $\CC^u_1$ and $\CC^v_2$ is equivalent to
the Codazzi equation for $(\hat A_{\hat\xi_2},\partial_u,\partial_v)$.
Notice that, as shown in the proof of \pref{tptqg},
Codazzi and Ricci equations for vectors in $\Delta$ and
$\Delta^\perp$ are just the projectability onto $L^2$ of the
functions involved and the determination of the splitting tensor,
 so they provide no additional
information.

For the last invariant, observe first that
$\tilde{\RR}=\cot(\theta)d\psi^1(\partial_{u},\partial_{v})$
is intrinsic. Indeed, by \eqref{chbag} and the Ricci equation,
$$
\tilde{\RR}=\cot(\theta)
\langle[A_{\xi_1},A_{\eta_1}]\partial_u,\partial_v\rangle=
\frac{\cos(\theta)}{\sin^2(\theta)}
\langle A_{\xi_1}\partial_v,A_{\xi_2}\partial_u\rangle=
\|\partial_u\|\|\partial_v\|\cos(\omega)\,{\rm scal},
$$
where for the last equality we used \eqref{gaussg}.
We conclude the invariance of $\RR$ from this and the previous
invariants, since it is straightforward to check that
$2\RR=2\tilde{\RR}+
(\CC^u_2)_v+(\CC^v_2)_u-(\ln {\cal G}^2)_{uv}$.
This completes the proof of the direct statement.

Conversely, given 5 smooth functions $\tau^u>0,\tau^v>0$,
$\hat\theta\in(0,\pi)$, $\hat\Lambda^u$ and $\hat\Lambda^v$ on $L^2$
satisfying \eqref{isg} for
$\hat\kappa^u=\tau^u\kappa^u$ and $\hat \kappa^v = \tau^v \kappa^v$, we
construct the second fundamental form and the normal connection 1-form
for a new isometric immersion as follows. We endow $T^\perp_fM$ with the
new metric that keeps $\hat\xi_1=\xi_1$ and $\hat\xi_2=\xi_2$ unitary
but making an angle~$\hat\theta$ instead of $\theta$. Then we set
\begin{equation}\label{As}
\hat A_{\hat\xi_1}=\sqrt{\tau^u}\ 
\frac{\sin(\hat\theta)}{\sin(\theta)}A_{\xi_1},\ \ \ \
\hat A_{\hat\xi_2}=\sqrt{\tau^v}\ 
\frac{\sin(\hat\theta)}{\sin(\theta)}A_{\xi_2},
\end{equation}
and
$$
\hat\psi^1(\partial_u):=
\tan(\hat\theta)\,\hat\Lambda^u-\hat\theta_u,\ \ \ 
\hat\psi^1(\partial_v)=-\tan(\hat\theta)\,\hat\Lambda^v,\ \ \
\hat\psi^1|_{\Gamma}=0.
$$
Therefore, as we saw in the direct statement, system \eqref{isg}, together
with the projectability of
$\tau^u,\,\tau^v$, $\hat \theta$, $\hat\Lambda^u$ and $\hat\Lambda^v$,
is equivalent to the fundamental equations for such a second fundamental
form and normal connection. Since $M^n$ is simply-connected, the
Fundamental Theorem of Submanifolds assures the existence of an isometric
immersion $\hat f:M^n\to\R^{n+2}$ with second fundamental form and
normal connection as above. It is clear by definition that
$\D(\hat f)=
\{\hat\theta,\hat\Lambda^u,\hat\Lambda^v, \hat\kappa^u,\hat\kappa^v\}$.
\qed

\section{The moduli space of deformations}\label{mmm}  

In this section we finally compute the moduli space of deformations of a
nowhere flat hyperbolic Euclidean submanifold $f$ by using the six
invariants found in the last section. Recall that, by \pref{diffnul} and
\pref{n12}, if such a deformation is not hyperbolic on some open subset,
then it is somewhere a composition and, in particular, it is not a honest
deformation.

\medskip

Let $f:M^n\to\R^{n+2}$ be a hyperbolic submanifold with main
angle $\theta$ and main symbols $\Lambda^u$ and $\Lambda^v$, as in
\eqref{lambdasdefs}. Fix $p_0=(u_0,v_0)\in L^2$, and fix primitives of
$\Lambda^u$ and $\Lambda^v$ with respect to $u$ and $v$ respectively,
that we denote as $\int\Lambda^udu:=\int_{u_0}^u\Lambda^u(t,v)dt
+\ln(s(u_0,v))/2$, and
$\int\Lambda^vdv:=\int_{v_0}^v\Lambda^u(u,t)dt+\ln(s(u,v_0))/2$.
In terms of a polar surface of $f$ we get from \eqref{lambdas} that
$\int\Lambda^udu=\ln(s)/2+\int_{u_0}^u(FG^{-1}\Gamma^u)(t,v)dt$ and
$\int\Lambda^vdv=\ln(s)/2+\int_{v_0}^v(FE^{-1}\Gamma^v)(u,t)dt$.

\medskip

Let $U=U(u)$ and $V=V(v)$ be a pair of functions of a real variable such
that
\begin{equation}\label{cfuvg}
U>-e^{-2\int\Lambda^vdv},\ V>-e^{-2\int\Lambda^udu},\ 
(1+Ue^{2\int\Lambda^vdv})(1+Ve^{2\int\Lambda^udu})>\cos^2(\theta).
\end{equation}
Using these two functions define $\tau^u$ and $\tau^v$ over $L^2$ as
\begin{equation}\label{tausg}
\tau^u:=1+Ve^{2\int\Lambda^udu}>0,\ \ \ \ \ \
\tau^v:=1+Ue^{2\int\Lambda^vdv}>0.
\end{equation}
Finally, define the second-order differential operator
$H_{UV}:C^\infty(L^2)\to C^\infty(L^2)$ by
\begin{equation}\label{H}
H_{UV}(\rho)=\rho_{uv}+\frac{\cos^2(\theta)}{\tau^u\tau^v}\,\rho_u\,\rho_v
-\frac{\Lambda^v}{\tau^v}\rho_u-\frac{\Lambda^u}{\tau^u}\rho_v.
\end{equation}

Our main result can be stated as follows. Observe that, as a consequence,
the moduli space of hyperbolic deformations of $f$ only depends on
$\Lambda^u,\Lambda^v$ and $\theta$, that are functions over the nullity
leaf space $L^2=M^n/\Gamma$. In particular, this moduli space only
depends on the metric of the polar surface~of~$f$.

\begin{theorem}\po\label{mainthg}
Let $f:M^n\to\R^{n+2}$ be a nowhere flat hyperbolic submanifold. Then,
the moduli space of local hyperbolic \isims of $M^n$ in
codimension two can be represented~as
$$
\moduli=\{(U,V):\eqref{cfuvg}\mbox{\normalfont{ holds, and }}
H_{UV}(\rho_{UV})=H_{00}(\rho_{00})\},
$$
where 
$\rho_{UV}
=\frac{1}{2}\ln\left(\frac{\tau^u\tau^v}{\cos^2(\theta)}-1\right)$.
\end{theorem}
\proof
Given a hyperbolic deformation $\hat f:M^n\to\R^{n+2}$,
we have by the invariance of
$\CC^u_1$, $\CC^v_1$ in \pref{ipg} that
$\hat\Lambda^u=\Lambda^u/\tau^u$ and $\hat\Lambda^v=\Lambda^v/\tau^v$.
This together with the invariance of $\CC^u_2$ and $\CC^v_2$
yields that $\tau^u$ and $\tau^v$ satisfy the uncoupled PDEs
\begin{equation}\label{pdetaug}
\tau^u_u=2\Lambda^u(\tau^u-1),\ \ \ \ \tau^v_v=2\Lambda^v(\tau^v-1),
\end{equation}
which are equivalent to the
existence of a pair of functions of one real variable $U(u)$ and $V(v)$
satisfying \eqref{tausg}. Furthermore, we have by the invariance of
${\cal G}$ that
\begin{equation}\label{cosdef}
\cos(\hat\theta)=\cos(\theta)/\sqrt{\tau^u\tau^v},
\end{equation}
\eqref{cfuvg} holds, and $\hat\rho:=\ln(|\tan(\hat\theta)|)=\rho_{UV}$
as in the statement. Thus, the invariance of $\RR$ is
equivalent to $H_{UV}(\rho_{UV})=H_{00}(\rho_{00})$, since
$f$ itself corresponds to $U=V=0$.

Conversely, if $(U(u),V(v))\in \moduli$, then \eqref{pdetaug} clearly
holds for $\tau^u>0$, $\tau^v>0$ given by \eqref{tausg} and consequently
$\hat\kappa^u=\tau^u\kappa^u$, $\hat\kappa^v=\tau^v\kappa^v$,
$\hat\Lambda^u=\Lambda^u/\tau^u$, $\hat\Lambda^v=\Lambda^v/\tau^v$
satisfy
$\CC^u_1=\hat\kappa^u\hat\Lambda^u$,
$\CC^v_1=\hat\kappa^v\hat\Lambda^v$,
$\CC^u_2=\hat\kappa^u_u/\hat\kappa^u+2\hat\Lambda^u$, and
$\CC^v_2=\hat\kappa^v_v/\hat\kappa^v+2\hat\Lambda^v$.
Moreover, $\hat\theta=\arccos(\cos(\theta)/\sqrt{\tau^u\tau^v})$
automatically satisfies that 
${\cal G}=\cos(\hat\theta)\sqrt{\hat\kappa^u\hat\kappa^v}$.
Lastly, setting $\hat\rho=\ln(|\tan(\hat\theta)|)=\rho_{UV}$, the
assumption $H_{UV}(\rho_{UV})=H_{00}(\rho_{00})$ means precisely that
$\RR=\hat\rho_{uv}+\cos^2(\hat\theta)\,\hat\rho_u\,\hat\rho_v
-\hat\Lambda^v\hat\rho_u-\hat\Lambda^u\hat\rho_v$.
Therefore, we conclude from \pref{ipg} that there exists locally
an isometric immersion $\hat f:M^n\to\R^{n+2}$ with
$\D(\hat f)=
\{\hat\theta,\hat\Lambda^u,\hat\Lambda^v,\hat\kappa^u,\hat\kappa^v\}$.
\qed

\begin{remark}\po\label{other}
{\rm Depending on the problem, it may be useful to change the operator in
\eqref{H} by composing $\rho$ with a suitable function. For example,
using the linear operator $\tilde H_{UV}(\rho)=
\rho_{uv}-\frac{\Lambda^v}{\tau^v}\rho_u-\frac{\Lambda^u}{\tau^u}\rho_v$,
we get in terms of
$2\tilde\rho_{UV}=\ln(1+\sin(\hat\theta))-\ln(1-\sin(\hat\theta))$ that
$\moduli=\{(U,V):\eqref{cfuvg}\mbox{\normalfont{ holds, and }}
\tilde H_{UV}(\tilde\rho_{UV})=
\tanh(\tilde\rho_{UV})\tilde H_{00}(\tilde\rho_{00})/2\sin(\theta)\}.$}
\end{remark}

\section{Honest and genuine deformations}\label{gensec} 

After classifying in the previous section the substantial rank two
deformations of a hyperbolic Euclidean submanifold in codimension two, we
now determine when such a deformation is a composition, when it is
genuine, and therefore when it is honest.

\medskip

\medskip

To describe the type of deformations of a rank two \isim we first need
the following result.

\begin{lemma}\po\label{gens}
Let $f,\hat f:M^n\to\R^{n+2}$ be a pair of nowhere congruent and nowhere
flat rank two \isims neither of which is contained in an affine
hyperplane. Then, there is an open dense subset $W\subset M^n$ along
which $f$ and $\hat f$ have two dimensional first normal spaces, and
either one of the following possibilities occur on each connected
component of $W$:
\begin{enumerate}[ 1)]
\item There are no unit normal vector fields $\mu$ of $f$ and $\hat \mu$
of $\hat f$ whose respective shape operators coincide. In this case,
the pair $\{f,\hat f\}$ is genuine;
\item There are orthonormal normal frames $\{\mu,\beta\}$ of $f$ and
$\{\hat \mu,\hat \beta\}$ of $\hat f$ such that $A_\mu=\hat A_{\hat\mu}$
and $\rk A_\beta = \rk \hat A_{\hat\beta} = k$, with $1\leq k\leq 2$.
In this case, $\{f,\hat f\}$ extends isometrically as
(regular) Sbrana-Cartan hypersurfaces if $k=2$, or as either flat or
(singular) generalized Sbrana-Cartan hypersurfaces if $k=1$.
In addition, when \mbox{$k=1$} and $f$ and $\hat f$ are not mutually
ruled, they extend isometrically as flat hypersurfaces if and only if
$\ker A_\beta \subset \ker\psi$, where $\psi = \la\nabla^\perp_\bullet
\mu,\beta\ra$ is the normal connection form of $\{\mu,\beta\}$.
\end{enumerate}
\end{lemma}
\proof
The dimension property of the first normal spaces is a consequence of
\pref{n12} and the discussion before it. In particular, if the second
fundamental forms of $f$ and $\hat f$ coincide in an open subset, then
by the Codazzi equation also their normal connections agree. Hence, the
immersions would be congruent along any simply-connected open subset
in~$U$. Thus, the second fundamental forms are almost everywhere
different.

Case (1) is obvious, since a normal vector field to the submanifold
tangent to the isometric extension gives a normal direction where the
shape operators coincide. So assume that there are unit normal vector
fields $\mu$ and $\hat \mu$ such that $A_\mu = \hat A_{\hat\mu}$, and
complete them to orthonormal normal frames $\{\mu,\beta\}$ of $f$ and
$\{\hat \mu,\hat \beta\}$ of $\hat f$. In particular, $1\leq\rk
A_\beta=\rk \hat A_{\hat\beta}\leq 2$ by the Gauss equation. The
smoothness of this frame is assured by Lemma 7 in \cite{dm2}.


If $\rk A_\beta=\rk \hat A_{\hat\beta}=2$, Lemma 6 in \cite{dm2} (or
Proposition 9 in \cite{df2}) says that case~(2) holds with regular
Sbrana-Cartan hypersurfaces as extensions.

Finally, assume that $\rk A_\beta=\rk \hat A_{\hat\beta}=1$. Lemma 9 in
\cite{dm2} assures that the pair extends as generalized Sbrana-Cartan
hypersurfaces, unless $\ker A_\beta \subset \ker\psi$, which is
equivalent for $A_\mu=\hat A_{\hat\mu}$ to be a Codazzi tensor by the
Codazzi equation for $A_\mu$. Hence, $\ker A_\beta \subset \ker\psi$ is
also a necessary condition for $f$ to extend as a flat hypersurface. In
particular, $\ker\hat A_{\hat\beta} \subset \ker\hat\psi$ also. We claim
that, in this situation, $f$ and $\hat f$ extend isometrically as flat
hypersurfaces.

To prove the claim, we have to consider two cases:

$(i)$ $A_{\beta}$ and $\hat A_{\hat\beta}$ are pointwise linearly
dependent. Define the pair $({\cal T},D)$, where $D:=\ker A_{\beta}$, and
${\cal T}$ is the line bundle isometry that sends $\mu$ into $\hat \mu$.
We conclude from Proposition 9 in \cite{df2} applied to $({\cal T},D)$
that $f$ and $\hat f$ extend isometrically as hypersurfaces with common
relative nullity of dimension $n=1+\dim D$, hence flat.

$(ii)$ $A_{\beta}$ and $\hat A_{\hat\beta}$ are pointwise linearly
independent. Notice that this implies also that
$\ker A_\beta \subset \ker\psi$ by comparing the Codazzi equations for
$A_\mu=\hat A_{\hat\mu}$. Now, we proceed as in $(i)$ defining
$({\cal T},D)$, but now for
$D:=\ker A_{\beta}\cap\ker\hat A_{\hat\beta}=\Delta$. By dimension
reasons, ${\cal N}(\phi)\supsetneq D$ for the bilinear form $\phi$
defined in equation (3) in \cite{df2}, and again by Proposition~9 in
\cite{df2} $f$ and $\hat f$ extend isometrically by relative nullity.
Thus, their extensions have rank one, hence flat.
\qed
\vspace{.5ex}

\medskip

Notice that, when $k=1$ in case (2) above, the condition
$\ker A_\beta \subset \ker\psi$ is equivalent for $f$ to be a composition.
In fact, if $e$ is an eigenvector field of $A_\beta$ with $A_\beta
e=\lambda e\neq 0$, then $\psi = \gamma \la e,\cdot\ra$ for some function
$\gamma$ and it is easy to check that $F(t,x)=f(x)+t(\gamma\lambda^{-1}
e+\mu)(x)$ is an immersion with relative nullity distribution
$\{\partial_t\}\oplus\ker A_\beta$. In particular $F$ has rank one, so it
is flat.

\begin{remark}\po\label{flatext}
{\rm By \pref{diffnul}, \pref{n12} and the above, if $M^n$ is a
rank two hyperbolic Riemannian manifold, then an \isim $f:M^n\to\R^{n+2}$
is a composition if and only if either $\dim N_f^1\neq 2$, or $\nu_f\neq
n-2$, or $\ker A_{\xi_i} \subset \ker\psi^i$, for some $i=1,2$. Observe
that, by \eqref{ker1g}, \eqref{2v1u} and \pref{cspg}, in terms of a polar
surface of $f$ the latter is equivalent to either $\Gamma^u=0$ or
$\Gamma^v=0$.}
\end{remark}

Let $f:M^n\to\R^{n+2}$ be a nowhere flat hyperbolic submanifold with main
angle $\theta$ and main symbols $\Lambda_u,\Lambda_v$ as in
\eqref{lambdasdefs} (or as in \eqref{lambdas} in terms of a polar
surface $g$ of $f$). In view of \rref{flatext} and \eqref{nc1fg}, the
condition for $f$ to be a composition is that either
$2s\Lambda^u=s_u$, or $2s\Lambda^v=s_v$. In particular, a hyperbolic
deformation $\hat f$ of $f$ is a composition if and only if either
$2\hat\Lambda^u=\ln(\sin^2(\hat\theta))_u$ or
$2\hat\Lambda^v=\ln(\sin^2(\hat\theta))_v$. But
\begin{eqnarray*}
2\hat\Lambda^u-\ln(\sin^2(\hat\theta))_u
&=&\frac{2\Lambda^u}{\tau^u}-
\ln\left(1-\frac{1-s}{\tau^u\tau^v}\right)_u\\
&=&-\frac{1}{\tau^u\tau^v-1+s}
\left(s_u-2\Lambda^u(\tau^v-1+s)+(1-s)\frac{\tau^v_u}{\tau^v}\right),
\end{eqnarray*}
and similarly for $2\hat\Lambda^v-\ln(\sin^2(\hat\theta))_v$.
By \eqref{pdetaug} we conclude that $\hat f$ is a
composition if and only if
\begin{equation}\label{notcomp}
\left(\frac{\tau^v-1+s}{\tau^v\,e^{2\int\Lambda^udu}}\right)_u=0
\ \ \ \ \ \ \ {\rm or}\ \ \ \ \ \ \ 
\left(\frac{\tau^u-1+s}{\tau^u\,e^{2\int\Lambda^vdv}}\right)_v=0.
\end{equation}

\medskip

We have the following criteria to extend isometrically both immersions
as flat hypersurfaces.

\begin{proposition}\po\label{extflathip}
Let $f:M^n\to\R^{n+2}$ be a nowhere flat hyperbolic submanifold and
$\hat f$ a hyperbolic deformation of $f$ given by $(U,V)\in\moduli$.
Then, in terms of a polar surface of $f$, the pair $\{f,\hat f\}$ extends
isometrically as flat hypersurfaces if and only if either $\Gamma^u=U=0$,
or $\Gamma^v=V=0$, or $\Gamma^u=V+1=0$, or $\Gamma^v=U+1=0$.
The latter two cases are equivalent to $\tau^u=\cos^2(\theta)$ and
$\tau^v=\cos^2(\theta)$, respectively.
\end{proposition}
\proof
According to \rref{flatext} and \eqref{notcomp}, a necessary condition
to extend as flat hypersurfaces is that either
$\Gamma^u\!=\!0$ or $\Gamma^v\!=\!0$, and either
$\left(\frac{\tau^v-1+s}{\tau^v\,e^{2\int\Lambda^udu}}\right)_u=0$ or 
$\left(\frac{\tau^u-1+s}{\tau^u\,e^{2\int\Lambda^vdv}}\right)_v=0$.
In order for the extensions to be isometric, we also need that the
corresponding second fundamental forms agree up to sign, as in
\lref{gens} part (2). So we have two
possibilities, up to obvious index choices:

$(i) \ A_{\eta_1}=\pm\hat A_{\hat \eta_1}$,
$\Gamma^u\!=\!0$, $\left(\frac{\tau^v-1+s}{s\tau^v}\right)_u\!=0$.
Here, from \eqref{chbag} we easily obtain that
the first equation is equivalent to $\tau^v=1$, that is, $U=0$, which
implies the third equation.

$(ii)\ A_{\eta_1}=\pm\hat A_{\hat \eta_2}$, $\Gamma^u\!=\!0$,
$\left(\frac{\tau^u-1+s}{\tau^u\,e^{2\int\Lambda^vdv}}\right)_v\!=0$.
In this case, the first equation is equivalent to $\tau^u=\cos^2(\theta)$
and then, by the second, $V=-1$ and so the third holds.

In any case, $\Gamma^u=0$, $\ker A_{\xi_i} \subset \ker\psi^i$ by
\rref{flatext}, and the proposition follows from \lref{gens} part (2)
since $f$ is not parabolic.
\qed
\vspace{1.5ex}

Our next principal result describes which deformations are genuine 
and honest, and how the ones that are not extend:

\begin{theorem}\po\label{genus}
Consider two nowhere congruent nowhere flat hyperbolic
\isims $f,\hat f:M^n\to\R^{n+2}$ which
do not extend isometrically as flat hypersurfaces.
Then, $\hat f$ is locally determined by a pair
$(U,V)\in \moduli$, and it holds that:
\begin{enumerate}[ 1)]
\item If $UV>0$, then $\{f,\hat f\}$ is genuine, and $\hat f$ is honest
if in addition \eqref{notcomp} does not hold;
\item If $UV=0$, then the pair $\{f,\hat f\}$ extends isometrically in a
unique way, and they do so as generalized (singular)
Sbrana-Cartan hypersurfaces. In particular, $\{f,\hat f\}$ is genuine;
\item If $UV<0$, then the pair $\{f,\hat f\}$ extends isometrically in
precisely two different ways, and they do so as (regular) Sbrana-Cartan
hypersurfaces. In particular, $f$ and $\hat f$ are constructed as the
intersection of two pairs of isometric Sbrana-Cartan hypersurfaces.
\end{enumerate}
Moreover, in cases (2) and (3), all the Sbrana-Cartan
extensions are of continuous or discrete class.
\end{theorem}
\proof
Since the immersions are hyperbolic, neither is contained in an affine
hyperplane. Moreover, as we already saw, the splitting tensor is
intrinsic and we have that $\hat A_{\hat\xi_i}$ and~$A_{\xi_i}$ are
linearly dependent for $i=1,2$. According to \lref{gens}, we have to
analyze when there are normal directions $\mu$ and $\hat \mu$ with the
same norm for which the corresponding shape operators coincide. Assume
this is the case, i.e.,
\begin{equation}\label{eq}
A_\mu=\hat A_{\hat\mu},
\end{equation}
and set
$\mu = a_1\xi_1+a_2\xi_2\neq 0$,
$\hat\mu = \hat a_1\hat\xi_1+\hat a_2\hat\xi_2\neq 0$, with
\begin{equation}\label{eqn}
a_1^2+a_2^2+2a_1a_2\cos(\theta)
=\hat a_1^2+\hat a_2^2+2\hat a_1\hat a_2\cos(\hat\theta).
\end{equation}
Evaluating \eqref{eq} in $0\neq Z_i\in\ker A_{\xi_i}\cap\Delta_f^\perp$,
we have
\begin{equation}\label{eai}
a_i\lambda_i=\hat a_i\hat \lambda_i, \ \ i=1,2.
\end{equation}
First observe that, if $\lambda_1=\hat \lambda_1$ and
$\lambda_2=\hat \lambda_2$, from \eqref{gaussg} we conclude that
$\theta=\hat \theta$, and thus $f$ and $\hat f$ would be congruent.
Hence assume that, say, $\lambda_2\neq\hat \lambda_2$, which also implies
that $a_1\neq 0$ and $\hat a_1\neq 0$ in view of \eqref{eai}. Now,
dividing \eqref{eqn} by $\hat a_1^2$ and using
\eqref{eai} we get for $t:=-\hat a_2/\hat a_1$ and
$k:=\sin^2(\hat\theta)/\sin^2(\theta)$ that
$$
t^2(k\tau^v-1)-2t(k\cos(\theta)\sqrt{\tau^u\tau^v}-\cos(\hat\theta))
+(k\tau^u-1)=0.
$$
This is a second degree polynomial since $\lambda_2\neq\hat \lambda_2$.
In view of \eqref{cosdef} this is equivalent to
\begin{equation}\label{poli}
h(t):=t^2(k\tau^v-1)-2t\cos(\hat\theta)(k\tau^u\tau^v-1) +(k\tau^u-1)=0,
\end{equation}
whose discriminant with respect to $t$ is
$-4k(\tau^u-1)(\tau^v-1)$,
which by \eqref{tausg} has the same sign as $-UV$.
Case (1) is then a consequence of case (1) in \lref{gens}.

For case (2), assume that $U=0$, and hence $\tau^v=1$.
Then $t=\cos(\hat\theta)^{-1}$ is the only root of \eqref{poli}.
Therefore, $\eta_1$ and $\hat\eta_1$ are the
only directions for which the shape operators coincide. But these are
precisely the directions orthogonal to $\xi_1$ and $\hat\xi_1$, whose
shape operators have rank one. We conclude case (2) from case (2) in
\lref{gens} for $k=1$.

For case (3), we have two different roots in \eqref{poli}, and we
claim that neither is equal to $\cos(\hat\theta)^{-1}$ or
$\cos(\hat\theta)$. To prove this, first we easily check that
\begin{equation}\label{cr}
h(\cos(\hat\theta)^{-1})=\frac{k\tau^u(\tau^v-1)}
{\cos^2(\theta)}(\tau^v-\cos^2(\theta)),\ \
h(\cos(\hat\theta))=\frac{k(\tau^u-1)}
{\tau^u}(\tau^u-\cos^2(\theta)).
\end{equation}
Since $UV\neq0$, then $\tau^u,\tau^v\neq 1$. But if, say,
$\tau^v=\cos^2(\theta)$, in view of \eqref{lambdas} and \eqref{tausg}
the function $e^{-2\int FE^{-1}\Gamma^v dv}=-U$ does not depend on $v$.
Since $M^n$ is nowhere flat, $F\neq 0$ and therefore $\Gamma^v=0$
and $U=-1$.
By \pref{extflathip} both $f$ and $\hat f$ extend isometrically as flat
hypersurfaces contradicting our hypothesis, and the claim is proved. The
proof of case (3) now follows from the discussion for case (2) and case
(2) in \lref{gens} for $k=2$, since the claim is equivalent to the fact
that $\mu$ and $\hat \mu$ are not collinear with $\eta_i$ and
$\hat \eta_i$, $i=1,2$.

\medskip

Now we argue that all the Sbrana-Cartan hypersurfaces that appear are
always of continuous or discrete class, that is, the
`interesting' classes III and IV in Theorem 3 in~\cite{dft2}. First,
recall that the extensions are by relative nullity, and then the relative
nullity of the codimension two hyperbolic submanifold $f$ is contained in
the relative nullity of its extension $F$, which also has rank two. In
particular, the splitting tensor $\tilde C_T$ of the relative nullity of
$F$ for $T\in \Delta_f\subset\Delta_F$ is conjugate to $C_T$. Since there
is $T\in\Delta_f$ such that $C_T$ has two different real eigenvalues, the
same holds for $\tilde C_T$, and hence, according to Theorem 3 in
\cite{dft2} for $c=0$, the extension is of continuous or discrete class.
\qed

\begin{remark}\po\label{2inter}{\rm
When $UV<0$ as in case (3), an interesting and unusual phenomenon occurs.
First, notice that, generically, intersections of rank two hypersurfaces
only provide rank 4 submanifolds. Yet, $f$ has rank two and is
constructed as the transversal intersection of a pair of non-isometric
Sbrana-Cartan hypersurfaces $L_1^{n+1},L_2^{n+1}\subset\R^{n+2}$, while
$\hat f$ is the transversal intersection of their respective deformations
\mbox{$\hat L_1^{n+1},\hat L_2^{n+1}\subset\R^{n+2}$}. That is,
\begin{equation}\label{exsbca}
f(M^n)=L_1^{n+1}\cap L_2^{n+1}\ \ \ \ \ \ {\rm and} \ \ \ \ \ \
\hat f(M^n)=\hat L_1^{n+1}\cap \hat L_2^{n+1}.
\end{equation}
In particular, although not a honest deformation, this provides
examples of interesting Sbrana-Cartan hypersurfaces of the continuous
or discrete classes. The lesson we extract from this is not to disregard
non genuine deformations, but instead use them to study deformability
in lower dimensions. We will see more examples of this kind of phenomena
in the last two sections.
}\end{remark}

\begin{remark}\po\label{comps}
{\rm As shown in its proof, when the pair $\{f,\hat f\}$ extends
isometrically as flat hypersurfaces, \tref{genus} still holds except in
the following two situations:
\begin{itemize}
\item In case (2), they may extend isometrically in a unique way but as
flat hypersurfaces instead of singular Sbrana-Cartan hypersurfaces, yet
if and only if either \mbox{$U=\Gamma^u=0$} or $V=\Gamma^v=0$, as seen in
\pref{extflathip};
\item In case (3), they extend isometrically as flat hypersurfaces
only if either $\Gamma^v=0$ and $U=-1$, or $\Gamma^u=0$ and $V=-1$, which
correspond to $\cos(\hat\theta)^{-1}$ or $\cos(\hat\theta)$ to be
roots of \eqref{poli}, respectively. But both cannot be
roots simultaneously, since otherwise by \eqref{cr} we would have
$\tau^u=\tau^v=\cos^2(\theta)$, which contradicts the third condition in
\eqref{cfuvg}. We conclude that $\{f,\hat f\}$ extends
isometrically also as (regular) Sbrana-Cartan hypersurfaces.
In other words, \eqref{exsbca} still holds, but with one of the pairs
$L^{n+1}_i,\hat L^{n+1}_i$ being flat, for some $i=1,2$.
\end{itemize}}
\end{remark}

\medskip

As shown in the Examples in \cite{df3} page 207, the singular set
$\Sigma^n$ of a hyperbolic Sbrana-Cartan hypersurface
$F:N^{n+1}\to\R^{n+2}$ is always a deformable rank two hyperbolic
Euclidean submanifold in codimension two, and $N^{n+1}$ itself can be
recovered from $\Sigma^n$. On the other hand, as a consequence of
\tref{genus} and \rref{comps}, we have that no pair $\{f,\hat f\}$ can
extend simultaneously both singularly and regularly. This answers
positively the natural question that was left open in \cite{dm2}, namely,
whether it is actually necessary to consider singular extensions to
obtain global genuine rigidity. Indeed, if~$F'$ is a deformation of $F$,
for any compact hypersurface $M^n\subset N^{n+1}$, we conclude that
$F(M^n)\subset\R^{n+2}$ and $F'(M^n)\subset\R^{n+2}$ are nowhere
congruent, yet they can only extend singularly along the interior of
$M^n\cap\Sigma^n$. Therefore, we have:

\begin{corollary}\po\label{neces}
The global genuine rigidity for compact Euclidean submanifolds in
codimension two as established in \cite{dm2} does not hold without
considering singular extensions.
\end{corollary}

\section{Hyperbolic submanifolds as hypersurfaces}\label{hsahg}  

In the recent paper \cite{dft3} the moduli space of all (local) isometric
immersions $f:M^n\to\R^{n+2}$ of a given Euclidean hypersurface
$g:M^n\to\R^{n+1}$ that are not compositions was computed. We can use the
machinery built in this work to understand the converse problem: to
classify rank two Euclidean submanifolds in codimension two that are also
hypersurfaces, and actually classify all their deformations. We will
carry out the study for hyperbolic submanifolds since these are
the ones that interest us in this paper, but, as we pointed out in the
introduction, similar analysis holds for the elliptic ones just by
taking complex conjugate coordinates instead of real ones, as done in
\cite{dft2} and \cite{dft3}.

\medskip

Let $f:M^n\to\R^{n+2}$ be a simply-connected nowhere flat hyperbolic
submanifold. In order to find an isometric immersion of $M^n$ as a
Euclidean hypersurface we will use the Fundamental Theorem of
Submanifolds by constructing a self-adjoint endomorphism $A$ on $TM$ that
satisfies the Gauss and Codazzi equations. Since $M^n$ has rank two, so
does~$A$, and $\Delta_f=\ker A$. Since $A_{\xi_1},A_{\xi_2}$ form a basis
of the self-adjoint tensors that satisfy \eqref{intcod}, we have that
$A=a_1 A_{\xi_1}+a_2 A_{\xi_2}$, where we can assume that $a_1>0$. By
\eqref{chbag} and \eqref{gaussg}, the Gauss equation for $A$ reduces to
$a_1a_2=-\cos(\theta)/\sin^2(\theta)$. So defining $\mu=a_1^2$ we have
$$
A=\sqrt{\mu}\, A_{\xi_1}-\frac{\cos(\theta)}{\sin^2(\theta)\sqrt{\mu}}
\, A_{\xi_2}.
$$
Using the notation $DB(X,Y):=\nabla_XBY-\nabla_YBX-B[X,Y]$ for (1,1)
tensors, and $\hat w(X,Y) := w(X)Y-w(Y)X$ for 1-forms, we have that
the Codazzi equation for $A$ is simply $DA=0$. So,
$$
\frac{1}{2}\, A_{\xi_1}\hat{d\mu} + {\mu}\, DA_{\xi_1}
+A_{\xi_2}\left(\frac{\cos(\theta)}{2\sin^2(\theta){\mu}}\hat{d\mu}+
\frac{1+\cos^2(\theta)}{\sin^3(\theta)}\hat{d\theta}\right)
-\frac{\cos(\theta)}{\sin^2(\theta)}\, DA_{\xi_2}=0.
$$
Recall that the Codazzi equation for $A_{\xi_i}$ is
$\sin(\theta)DA_{\xi_i}=(-1)^j(A_{\xi_j}-\cos(\theta)A_{\xi_i})\hat\psi^i$,
for $1\leq i\neq j \leq2$. Hence, since the images of $A_{\xi_1}$ and
$A_{\xi_2}$ are linearly independent we get
$$
A_{\xi_1}\!\left(
\frac{1}{2}\hat{d\mu}-\mu\frac{\cos(\theta)}{\sin(\theta)}\hat\psi^1
+\frac{\cos(\theta)}{\sin^3(\theta)}\hat\psi^2\right)=0,
$$
$$
A_{\xi_2}\!\left(\frac{1}{2{\mu}}\hat{d\mu}
+\mu\frac{\sin(\theta)}{\cos(\theta)}\hat\psi^1
-\frac{\cos(\theta)}{\sin(\theta)}\hat\psi^2
+\frac{1+\cos^2(\theta)}{\cos(\theta)\sin(\theta)}\hat{d\theta}
\right)=0.
$$

Observe first that these equations for one vector in $\Delta$ and the
other in $\Delta^\perp$ say that~$\mu$ is projectable, since
$\psi^1,\psi^2$ and $\theta$ also are. In view of \eqref{nc1fg}, we
obtain that the above two equations are equivalent to the first order
system of PDE
\begin{equation}\label{hh}
\mu_u=\mu a-b,\ \ \mu_v=\mu(\mu c-d),
\end{equation}
where, as usual, $s=\sin^2(\theta)$, and
$$
a=2\Lambda^u-\frac{s_u}{s},\ \
b=\frac{2\Lambda^u}{s},\ \
c=\frac{2s\Lambda^v}{1-s},\ \
d=2\Lambda^v+\frac{s_v}{s(1-s)}.
$$
The integrability condition of \eqref{hh} is therefore
\begin{equation}\label{ic}
    P(\mu):=(ac+c_u)\mu^2-(2bc+d_u+a_v)\mu+(bd+b_v)=0.
\end{equation}
We point out for further reference that
\begin{eqnarray*}
ac+c_u&=&\frac{2s}{1-s}\left(\Lambda^v_u+2\Lambda^u\Lambda^v
+\Lambda^v\frac{s_u}{1-s}\right),\\
2bc+d_u+a_v&=&
2\Lambda^v_u+2\Lambda^u_v+\frac{1}{1-s}\left(8\Lambda^u\Lambda^v+s_{uv}
+\frac{s_us_v}{1-s}\right),\\
bd+b_v&=&\frac{2}{s}\, \left(\Lambda^u_v+2\Lambda^u\Lambda^v
+\Lambda^u\frac{s_v}{1-s}\right).
\end{eqnarray*}

We conclude that one and only one of the following possibilities,
enumerated from the least to the most generic, holds
along each connected component of an open dense subset of $M^n$:
\begin{enumerate}[$i)$]
\item $P=0$, that is, $a=-\ln(|c|)_u, d=-\ln(|b|)_v$ and
$2bc=\ln(|bc|)_{uv}$, in which case the manifold admits an
\isim as a Sbrana-Cartan hypersurface of the continuous class;
\item Equation \eqref{ic} has two positive roots, and both satisfy
\eqref{hh}, in which case the manifold admits an \isim
as a Sbrana-Cartan hypersurface of the discrete class;
\item Only one of the positive roots of \eqref{ic} satisfies \eqref{hh},
in which case the manifold admits an \isim as a rigid hypersurface;
\item No positive root of \eqref{ic} satisfies \eqref{hh}, hence the
manifold admits no \isim as a hypersurface.
In particular, this is the case if
$(2bc+d_u+a_v)^2<4(ac+c_u)(bd+b_v)$, or $2bc+d_u+a_v\leq0,\,bd+b_v\geq0,\,
ac+c_u\geq0$, or $2bc+d_u+a_v\geq0,\,bd+b_v\leq0,\,ac+c_u\leq0$.

\end{enumerate}

\section{Deformations preserving the main angle}\label{mangle}  

This and the following section are devoted to give some applications. The
purpose in this one is to describe a particularly interesting class of
deformations of a hyperbolic submanifold: the ones preserving the main
angle $\theta$. Although we will see that these deformations are never
honest nor genuine, they provide interesting applications to the
Sbrana-Cartan theory of deformable hypersurfaces. This justifies what we
pointed out in the introduction: one should not simply ignore the study
of non honest and non genuine deformations since they can provide
insights for lower codimension rigidity.

\bigskip

So, we study here the implications of \tref{mainthg} on the
structure of a hyperbolic nowhere flat submanifold $f:M^n\to\R^{n+2}$
admitting a hyperbolic deformation $\hat f$ with
$\hat\theta=\theta$. We know that $\hat f$ is determined by
$(U,V)\in \moduli$, and by \eqref{cosdef} the condition which
characterizes these deformations is simply that
\begin{equation}\label{gausspag}
\tau:=\tau^u=1/\tau^v.
\end{equation}
Observe that $UV<0$ since $(\tau^u-1)(\tau^v-1)=-(\tau-1)^2/\tau<0$, and
$\tau\neq 1$ since the immersions are not congruent. Hence, in view of
\rref{comps}, \tref{genus} holds and $f$ and $\hat f$ extend
isometrically as Sbrana-Cartan hypersurfaces in two different ways as in
\rref{2inter}, unless $\tau=\cos^2(\theta)$ or $\tau=\cos^{-2}(\theta)$,
in which case $f$ and $\hat f$ extend isometrically both as flat and
Sbrana-Cartan hypersurfaces in a unique way. In any case, $\hat f$
is never a genuine deformation and, in fact, $\pm\sqrt\tau$ are the two
real roots of \eqref{poli}. Moreover, using \eqref{gausspag} in
\eqref{tausg} we obtain that $V=(\tau-1)\,e^{-2\int\Lambda^udu}$, and
$U=(\tau^{-1}-1)\,e^{-2\int\Lambda^vdv}$ depend on one variable only.
Equivalently,
\begin{equation}\label{sbcasy}
\tau_u=2\Lambda^u(\tau-1),\ \ \ \tau_v=2\Lambda^v\tau(\tau-1).
\end{equation}

\begin{remark}\po\label{sbcapres}
{\rm Observe that system \eqref{sbcasy} is exactly the system that
appears in the Sbrana-Cartan theory, but now for the Euclidean polar
surface of $f$ instead of the spherical Gauss map of its extension. Its
integrability condition is also given by
\begin{equation}\label{icsc}
\tau(\Lambda^v_u +2\Lambda^u\Lambda^v)=\Lambda^u_v +2\Lambda^u\Lambda^v,
\end{equation}
as in the Sbrana-Cartan theory.}
\end{remark}

To compute $\moduli$ we have $\rho_{00}=\rho_{UV}=\ln(|\tan(\theta)|)$
and thus
$$
0=H_{UV}(\rho_{UV})-H_{00}(\rho_{00})=
(1-\tau)\Lambda^v(\rho_{00})_u+(1-\tau^{-1})\Lambda^u(\rho_{00})_v.
$$
Since $\tau\neq 1$ this is equivalent to
\begin{equation}\label{d1expg}
\tau\Lambda^vs_u=\Lambda^us_v.
\end{equation}

First, consider this equation under the generic condition
$\Lambda^us_v\neq0$. Hence we also have $\Lambda^vs_u\neq0$. We conclude
that, in this generic situation, $f$ admits at most one deformation
$\hat f$ preserving the main angle, depending on whether
$1\neq\tau=\Lambda^us_v/\Lambda^vs_u>0$ satisfies system
\eqref{sbcasy} or not.

\medskip

Let us now turn our attention to the non generic case where
\eqref{d1expg} trivially holds, i.e.,
$$
\Lambda^us_v=\Lambda^vs_u=0.
$$
Thus $H_{UV}(\rho_{UV})=H_{00}(\rho_{00})$ is automatically satisfied,
and three possibilities may occur:
\begin{enumerate}[$\ i)$]
\item Either \ $s_u=\Lambda^u=0$,\ \ or\ \ $s_v=\Lambda^v=0$;
\item $\Lambda^v=\Lambda^u=0$;
\item $\Lambda^u\neq 0$, $\Lambda^v\neq 0$, and $\theta$ is constant.
\end{enumerate}

\medskip

{\it Case} $(i)$. Suppose that, say, $s_u=\Lambda^u=0$. Then, by
\eqref{sbcasy}, $\tau$ and $\Lambda^v$ are functions of $v$ only,
$V=\tau-1$ and $U$ is constant. We conclude from \tref{mainthg} that, in
this situation, a deformation of $f$ preserving $\theta$ exists if and
only if either $s_u=\Lambda^u=\Lambda^v_u=0$, or
$s_v=\Lambda^v=\Lambda^u_v=0$, in which case there is actually a
one-parameter family of such deformations, one for each constant chosen
for $U$ or $V$, respectively.

Observe that, in this case, \eqref{nc1fg} also gives
$\psi^1(\partial_{u})=0$, and hence $\ker A_{\xi_1}\subset\ker\psi^1$, or
analogously $\ker A_{\xi_2}\subset\ker\psi^2$. This implies by
\rref{flatext} that $f$ is a composition. Actually, for such an $f$ the
polynomial $P$ in \eqref{ic} vanishes identically. Therefore, $M^n$
admits also a one-parameter family of isometric deformations as a
Euclidean hypersurface, so each such an $M^n$ provides an
example of a Sbrana-Cartan hypersurface of the continuous class.

\medskip

{\it Case} $(ii)$. Here \eqref{sbcasy} simply says that $1\neq \tau>0$
is constant, $U=\tau^{-1}-1, V=\tau-1$. We conclude that every member of
the class of hyperbolic submanifolds satisfying $\Lambda^v=\Lambda^u=0$
admits a one-parameter family of deformations preserving the main angle.

Moreover, in this case, $P$ in \eqref{ic} is 
$P(\mu)=-\mu(s_{uv}/(1-s)+s_us_v/(1-s)^2)$ and thus $M^n$ admits no
\isim as a Euclidean hypersurface unless $(s-1)s_{uv}=s_us_v$, in
which case $M^n$ is also an example of a Sbrana-Cartan hypersurface of
the continuous class.

\medskip

{\it Case} $(iii)$. Here, we search for a function $\tau$ satisfying
\eqref{sbcasy}, whose integrability condition is \eqref{icsc}.
So, we have two subcases:

If $\Lambda^v_u +2\Lambda^u\Lambda^v\neq 0$ and $\Lambda^u_v
+2\Lambda^u\Lambda^v\neq 0$ are different and have the same sign,
then~$f$ admits at most one deformation $\hat f$ preserving the main
angle, depending on whether
$1\neq\tau=(\Lambda^u_v +2\Lambda^u\Lambda^v)
/(\Lambda^v_u +2\Lambda^u\Lambda^v)>0$
satisfies system \eqref{sbcasy} or not.

If, on the contrary, we have 
$\Lambda^v_u=\Lambda^u_v=-2\Lambda^u\Lambda^v\neq 0$,
this easily implies that $\Lambda^v=\tilde V'/2(\tilde U+\tilde V)$ and
$\Lambda^u=\tilde U'/2(\tilde U+\tilde V)$ for some non-constant
one-variable functions $\tilde U=\tilde U(u)$ and $\tilde V=\tilde
V(v)$. Then, it is easy to check that the pairs $(U,V)$ are given by
$U=1/(c-\tilde U)$ and $V=-1/(c+\tilde V)$, for $c\in\R$. Therefore, a
one-parameter family of deformations preserving $\theta$ always exists
in this case. Observe that $P$ in \eqref{ic} is
$P(\mu)=8\mu\Lambda^u\Lambda^vs/(1-s)$, it has no positive roots, and so
$M^n$ is not a Euclidean hypersurface.

\medskip

Summarizing, we have shown:

\begin{theorem}\po\label{prest}
A nowhere flat hyperbolic submanifold $f:M\to\R^{n+2}$ has a hyperbolic
deformation $\hat f$ preserving the main angle $\theta$ if and only if
either one of the following occurs:
\begin{enumerate}[\ \ 1)]
\item $\Lambda^us_v\neq0,\Lambda^vs_u\neq0$ and the function
$1\neq\tau=\Lambda^us_v/\Lambda^vs_u>0$ satisfies \eqref{sbcasy};
\item $\theta$ is constant,
$\Lambda^v_u +2\Lambda^u\Lambda^v\neq 0$,
$\Lambda^u_v +2\Lambda^u\Lambda^v\neq 0$, and the function
$1\neq\tau>0$ given by
$\tau=(\Lambda^u_v+2\Lambda^u\Lambda^v)/
(\Lambda^v_u +2\Lambda^u\Lambda^v)$ satisfies \eqref{sbcasy};
\item Either $s_u=\Lambda^u=\Lambda^v_u=0$, or
$s_v=\Lambda^v=\Lambda^u_v=0$, or $\Lambda^v=\Lambda^u=0$, or $\theta$ is
constant and $\Lambda^v_u=\Lambda^u_v=-2\Lambda^u\Lambda^v$.
\end{enumerate}
Moreover, in cases (1) and (2) $f$ has only one noncongruent such
deformation, while in case (3) it has precisely a one-parameter family of
them.

In any case, $\{f, \hat f\}$ extend isometrically precisely in two
different ways as in \eqref{exsbca}, both as Sbrana-Cartan hypersurfaces
if $\tau\neq1-s,(1-s)^{-1}$, or as Sbrana-Cartan hypersurfaces and flat
hypersurfaces otherwise. All these Sbrana-Cartan hypersurface extensions
are of the continuous or discrete classes. \end{theorem}

\section{Sbrana-Cartan hypersurfaces of intersection type}\label{sbcait}  

The main results in \cite{dft2}, Theorems 9 and 11, were devoted to the
construction of a large family of Sbrana-Cartan hypersurfaces of the
discrete class in any dimension by intersecting two flat hypersurfaces in
general position $N_i^{n+1}\subset\R^{n+2}$, that is,
\begin{equation}\label{intfh}
M^n=N_1^{n+1}\cap N_2^{n+1}\subset\R^{n+2}.
\end{equation}
These Sbrana-Cartan hypersurfaces $M^n\subset N_i^{n+1}\subset\R^{n+1}$,
which we call here {\it of intersection type}, are characterized by the
fact that their Gauss map satisfies $\Gamma^1_u-\Gamma^1\Gamma^2+F=0$;
see Lemma 10 in \cite{dft2}. It is immediate that these, as submanifolds
in codimension two $M^n\subset\R^{n+2}$, are nowhere flat and hyperbolic.
We proceed now to easily recover these two main results by using the
machinery developed in this work. Moreover, we classify
all their deformations in codimension two, finding the first known
examples of honestly deformable submanifolds in codimension two.

\medskip

An equivalent way to understand the Sbrana-Cartan hypersurfaces of
intersection type is to consider an embedded nowhere flat hyperbolic
submanifold $M^n\subset\R^{n+2}$, and ask for it to extend as flat
hypersurfaces in two different ways. As we saw in \rref{flatext}, this is
equivalent for its polar surface $g$ to satisfy
$\Gamma^u=\Gamma^v=0$ in \eqref{guv}, that is,
$$
g(u,v)=\alpha_1(u)+\alpha_2(v)
$$
is the sum of two regular curves, with
$\alpha_1',\alpha_1'',\alpha_2',\alpha_2''$ pointwise linearly
independent. To avoid surface-like submanifolds, we require also for $g$
to have conformal substantial codimension at least 3. We can further
assume that $\alpha_1$ and $\alpha_2$ are parametrized by arc-length,
i.e., $E=G=1$. Of course, in this situation
$F=\cos(\theta)=\la\alpha_1',\alpha_2'\ra$ is the sum of $n+2$ arbitrary
functions whose logarithms separate variables. We also have by
\eqref{lambdas} that
$$
2\Lambda^u = s_u/s,\ \ 2\Lambda^v = s_v/s.
$$
Observe in addition that $\xi_i=\alpha_j'$, and
$(-1)^i\sin(\theta)\eta_i=\cos(\theta)\alpha_j'-\alpha_i'$,
$1\leq i\neq j\leq2$.

\bigskip

We proceed to recover the main results in \cite{dft2} for which we
use the concept of local shared dimension $I$ of a pair of curves defined
at the end of \sref{prelg}.

\begin{theorem}\po\label{sc}
Let $i:M^n\subset\R^{n+2}$ be a nowhere flat embedded hyperbolic
submanifold, and assume its polar surface $g$ separates variables, i.e.,
$g(u,v)=\alpha_1(u)+\alpha_2(v)$. Then, $M^n$ is the transversal
intersection of two flat hypersurfaces as in \eqref{intfh}.
Moreover, as a hypersurface, $M^n\subset \R^{n+1}$ is a Sbrana-Cartan
hypersurface of the discrete class if $I(i):=I(\alpha_1,\alpha_2)\geq 2$,
and of the continuous class if $I(i)=1$.

Conversely, the polar surface of a Sbrana-Cartan hypersurface of
intersection type in codimension two separates variables.
\end{theorem}
\proof
We have already argued for the converse statement. For the direct one, in
our situation, the polynomial $P$ in \eqref{ic} is
\begin{equation}\label{sb1}
P(\mu)=\frac{1}{1-s}\left(s_{uv}+\frac{s_us_v}{1-s}\right)
\left(\mu-\frac{1}{s}\right)\left(\mu+1-\frac{1}{s}\right).
\end{equation}
In particular, $P(\mu)=0$ for $\mu=1/s$ and $\mu=1/s-1$. In other
words, $A_{\eta_i}$ satisfies the Gauss and Codazzi equations for
Euclidean hypersurfaces, $i=1,2$. We have two possibilities:

$P=0$. This is the case when
$(1-s)^{-2}((1-s)s_{uv}+s_us_v)=-(\ln(1-s))_{uv}=0$, or, equivalently,
$F=\cos(\theta)=a(u)b(v)$ is the product of two functions of one
variable. Hence, by \lref{cos}, $I(i)=1$. Now, the discussion at the end
of \sref{hsahg} already implies that $M\subset \R^{n+1}$ is a
Sbrana-Cartan hypersurface of the continuous class. But here we can do
better and actually solve \eqref{hh}: $\mu$ is given by
$\mu=\lambda(v)+1/s$, where $\lambda$ is any solution of the ODE\ \
$\lambda'=\lambda(\lambda+1)b'/b$\ \ for which $\mu>0$. Notice also that
$\lambda=0$ and $\lambda=-1$ are two solutions of this ODE, which
correspond to the two original intersecting flat hypersurfaces.

$P\neq 0$. In this situation \eqref{sb1} has precisely the
two positive solutions just described, $\mu=1/s,1/s-1$, and hence
$M^n\subset N_i^{n+1}\subset\R^{n+1}$ are the two unique noncongruent
\isims of $M^n$ as a Euclidean hypersurface.
\qed
\vspace{1.5ex}

As another application, we now compute all the honest deformations
$\hat i:M^n\to\R^{n+2}$ of any Sbrana-Cartan hypersurface of intersection
type.

\medskip

Suppose there is such an \isim $\hat i$. Since it is not a composition
it has rank two and by \tref{mainthg} it is induced by
$(U,V)\in\modulii$. In this case, $\tau^u = 1+sV(v)$ and
$\tau^v=1+sU(u)$, and thus the condition \eqref{notcomp} for such an
$\hat i$ not to be a composition turns out to be
$$
U'\neq -U(U+1)\ln(\cos^2(\theta))_u\ \ \ \ \ {\rm and}\ \ \ \ \
V'\neq -V(V+1)\ln(\cos^2(\theta))_v.
$$
In particular, we assume that $U,V\neq 0,-1$. Moreover, conditions
\eqref{cfuvg} are simply
\begin{equation}\label{conduv}
U,V>-1/s\ \ \ \ {\rm and}\ \ \ \
(U+1)(V+1)>\cos^2(\theta)UV.
\end{equation}
Define
\begin{equation}\label{uvt}
\varphi=\cos^2(\theta)\tilde U\tilde V,\ \ \ {\rm for}\ \ \
\tilde U=\frac{U}{U+1},\ \  \tilde V=\frac{V}{V+1}.
\end{equation}
We have that $\hat i$ is not a composition if and
only if $\varphi_u\neq0$ and $\varphi_v\neq0$.
It is easy to check that the equation defining $\modulii$ in
\tref{mainthg} now becomes
\begin{equation}\label{sepvardefs}
2\varphi(1-\varphi)\varphi_{uv} +(2\varphi-1)\varphi_u\varphi_v=0.
\end{equation}
According to \tref{sc} and \lref{cos}, $M^n\subset \R^{n+1}$ is a
Sbrana-Cartan hypersurface of the continuous class if and only if
$$
0=\ln(\cos^2(\theta))_{uv}=\ln(\varphi)_{uv}=\varphi^{-2}(\varphi\varphi_{uv}-\varphi_u\varphi_v).
$$
But then \eqref{sepvardefs} reduces to
$\varphi_u\varphi_v=0$. We conclude that the only deformations in
codimension two of these Sbrana-Cartan hypersurfaces of the continuous
class, for which $I(i)=1$, are compositions.

\medskip

So let us concentrate on the discrete class, i.e.,
$I(i)\geq 2$. If we set
$\tilde\varphi=\arcsin(2\varphi-1)$ when $\varphi\in(0,1)$ and
$\tilde\varphi=\ln(|2\varphi-1+2\sqrt{\varphi(\varphi-1)}|)$ otherwise,
\eqref{sepvardefs} is just $\tilde\varphi_{uv}=0$.
We claim that $I(i) = 2$.
Indeed, if
$\varphi\in(0,1)$ there are functions $U_0(u),V_0(v)$ such that
$2\varphi=\sin(U_0+V_0)+1$, and then
$$
2\la\alpha_1',\alpha_2'\ra^2=2\cos^2(\theta)
=(\tilde U\tilde V)^{-1}(\sin(U_0+V_0)+1) =
(\tilde U_1\tilde V_1\pm\tilde U_2\tilde V_2)^2,
$$
for $\tilde U_i^2=(1+(-1)^i\sin(U_0))/|\tilde U|$,
$\tilde V_i^2=(1+(-1)^i\cos(V_0))/|\tilde V|$, $i=1,2$.
A similar computation holds for $\varphi\not\in(0,1)$, and the claim
follows from \lref{cos}. In particular, for $I(i)\neq 2$, $i$ is
honestly rigid.

Then, assume from now on that $I(i) = 2$ is constant and, for
$j=1,2$, denote by $\overline \alpha_j$ the orthogonal projection of
$\alpha_j$ to the shared plane $\V^2$ between $\alpha_1$ and $\alpha_2$.
First, we claim that \eqref{sepvardefs} holds for
$\tilde U=t\|\overline \alpha_1'\|^{-2}$,
$\tilde V=t^{-1}\|\overline \alpha_2'\|^{-2}$ and $0\neq t\in\R$. Indeed,
in this situation, $\varphi=\la e_1,e_2\ra^2=\cos^2(w)$, where
$w=\angle(e_1,e_2)$ and
$e_j=\overline \alpha_j'/\|\overline \alpha_j'\|\in\V^2\cong\C$,
$j=1,2$, lies in the unit circle $\Sp^1\subset \C$. Hence, writing
$e_j'=k_j i e_j$ for a function $k_j\neq 0$ since $I(i)= 2$, we
have
$$
2\varphi(1-\varphi)\varphi_{uv} +(2\varphi-1)\varphi_u\varphi_v=
4k_1k_2\cos^2(w)\sin^2(w)(\cos(2w)-2\varphi+1)=0,
$$
and the claim is proved. Notice, in particular, that
$(\tilde U\tilde V)^{-1}>\cos^2(\theta)$, since
\begin{equation}\label{cash}
\|\overline\alpha_1'\|^2\|\overline\alpha_2'\|^2 >
\la\overline\alpha_1',\overline\alpha_2'\ra^2 =
\la\alpha_1',\alpha_2'\ra^2=\cos^2(\theta).
\end{equation}
We prove next that these are in fact all the solutions.

To determine all $\tilde U,\tilde V$ in \eqref{uvt} that satisfy
\eqref{sepvardefs}, write
$\tilde U(u) = \epsilon_1e^{U_1(u)}\|\overline\alpha_1'(u)\|^{-2}$
and $\tilde V(v) = \epsilon_2e^{V_1(v)} \|\overline \alpha_2'(v)\|^{-2}$,
and thus $\varphi = \epsilon_1\epsilon_2e^{U_1+V_1}\cos^2(w)$,
where $\epsilon_1,\epsilon_2=\pm1$. Since \eqref{sepvardefs} is
independent under change of parametrizations in $u$ and in $v$, we
can assume that $k_1=k_2=1$ and so $w=u-v\neq 0$.
Thus \eqref{sepvardefs} is
\begin{equation}\label{n1}
U_1'V_1'+2\tan(u-v)(U_1'-V_1')+4(1-\epsilon_1\epsilon_2e^{U_1}e^{V_1})=0.
\end{equation}
We claim that $U_1$ and $V_1$ are constant, with $U_1+V_1=0$, and
$\epsilon_1=\epsilon_2$. To prove the claim, observe that we can write
\eqref{n1} as $U_1(u)'=a(u,v)e^{U_1(u)}-b(u,v)$, with
$a=4\epsilon_1\epsilon_2e^{V_1}/(V_1'+2\tan(u-v))$,
$b=(4-2\tan(u-v)V_1')/(V_1'+2\tan(u-v))$, and similarly for $V_1'$.
Then, $e^{U_1}=b_v/a_v$ does not depend on $v$, and similarly for $V_1$.
A straightforward computation shows that $(b_v/a_v)_v=0$ if and only if
$\tan(u)(A+B\tan(v)) + (B-A\tan(v))=0$,
where $A=V_1'''+V_1''V_1'$, and $B=V_1'''V_1'/2-V_1''(V_1''+2)$. This
happens only when $A=B=0$ or, equivalently, if either $V_1''=0$ or
$V_1''=-V_1'^2/2-2$. Similarly, $U_1''=0$ or $U_1''=-U_1'^2/2-2$.
It is now easy to verify that the only possibility is that $U_1=-V_1$
is a constant and $\epsilon_1\epsilon_2=1$, as wished.

This claim implies that the set of deformations $\hat i=i_t$ of $i$
that are not compositions is the one-parameter family
$$
\left(U_t=(t^{-1}\|\overline \alpha_1'\|^2-1)^{-1},\ 
V_t=(t\|\overline \alpha_2'\|^2-1)^{-1}\right)\in\moduli
$$
satisfying \eqref{conduv}. In view of \eqref{cash}, conditions
\eqref{conduv} are equivalent to
$$
\frac{1}{1-t^{-1}\|\overline \alpha_1'\|^2}<\frac{1}{s},\ \
\frac{1}{1-t\|\overline \alpha_2'\|^2}<\frac{1}{s},\ \ \ {\rm and}\ \ \
(1-t^{-1}\|\overline \alpha_1'\|^2)(1-t\|\overline \alpha_2'\|^2)>0.
$$
These are obviously satisfied for all $t<0$. For $t>0$, it cannot happen
simultaneously that $t^{-1}\|\overline \alpha_1'\|^2<1$
and $t\|\overline \alpha_2'\|^2<1$ since, by the second equation,
$\|\overline \alpha_1'\|^2<t<\cos^2(\theta)\|\overline \alpha_2'\|^{-2}$,
contradicting \eqref{cash}.
On the other hand, $t^{-1}\|\overline \alpha_1'\|^2>1$
and $t\|\overline \alpha_2'\|^2>1$ is not possible either because
$\|\overline \alpha_1'\|^2>t>\|\overline \alpha_2'\|^{-2}$ contradicts
the fact that $\|\overline \alpha_j'\|^2\leq\|\alpha_j'\|^2=1$, $j=1,2$.
We conclude that $t<0$, and therefore $UV>0$. That is, by \tref{genus} we
obtain that the moduli space of hyperbolic deformations is a connected
differentiable 1-parameter family of honest deformations, and
therefore $i$ is genuinely and honestly deformable.

It is interesting to analyze the boundary of this
family, i.e., $t=0$ and $t=\infty$. We have in this case that $U=0,V=-1$
and $U=-1,V=0$, respectively. These do not satisfy the third condition in
\eqref{cfuvg}, yet they give rise to a pair of rank two \isims
$i_0,i_\infty$ of $M^n$ that lie inside a hyperplane. Indeed, say for
$V=0,U=-1$, we have at the limit that $\tau^u=\cos^2(\theta)$,
$\tau^v=1$, $\hat \theta=0$, $\hat\xi_1=\hat\xi_2$, $\hat\psi^i=0$,
$\hat A_{\hat\xi_i}=0$ by \eqref{As}, and in view of \eqref{chbag},
$$
\hat A_{\hat \eta_i}=
\lim_{t\to 0}\frac{(-1)^i}{\sin(\hat\theta)}
\left(\cos(\hat\theta)\hat A_{\hat\xi_i}-\hat A_{\hat\xi_j}\right)=
\lim_{t\to 0}\frac{(-1)^i}{\sin(\theta)}
\left(\cos(\hat\theta)\sqrt{\tau^i}A_{\xi_i}-\sqrt{\tau^j}A_{\xi_j}\right)=
A_{\eta_1},
$$
that is a Codazzi tensor. Similarly for $U=0,V=-1$ we get
$\hat A_{\hat \eta_i}=A_{\eta_2}$. In other words, $i_0$ and $i_\infty$
are precisely the two unique \isims of $M^n$ as a Euclidean
hypersurface. In particular, the two pairs $\{i,i_0\}$ and
$\{i,i_\infty\}$ are not genuine since each pair extends isometrically
in a unique way, and as flat hypersurfaces.

\medskip

Summarizing, we have shown:

\begin{theorem}\po\label{scdef}
A Sbrana-Cartan hypersurface of intersection type $i:M^n\subset\R^{n+2}$
is honestly rigid, unless $I(i)=2$. In the latter case, the moduli space
of local rank two deformations of $i$ is a differentiable compact
connected 1-parameter family $\{i_t:-1\leq t\leq1\}$. Moreover, the
interior members of this family, $i_t$ for $-1<t<1$, are honest
deformations of $i$, while the pair of deformations $i_{\pm1}$ at its
boundary are the only two \isims of $M^n$ as a Euclidean hypersurface.
In addition, $\{i,i_{\pm1}\}$ extend isometrically as flat hypersurfaces,
and therefore $i_{\pm1}$ are not genuine deformations of $i$.
\end{theorem}

\begin{remark}\po
{\rm It was not known until now if a honestly locally deformable
Euclidean submanifold of rank two in codimension two existed at all,
since, to our surprise, even the highly degenerate elliptic and parabolic
Euclidean submanifolds in codimension two were shown to be honestly rigid
in Theorems 1 and 4 in \cite{df3}. Now \tref{scdef} answers affirmatively
this question.}
\end{remark}

\tref{scdef} also shows that different kinds of deformations can be glued
smoothly in complex ways. Indeed, by taking two curves
$\alpha_1,\alpha_2$ for which $I(\alpha_1,\alpha_2)$ varies from point
to point, we can construct a connected submanifold that has an open dense
subset such that each connected component deforms in different ways.

\medskip


{\renewcommand{\baselinestretch}{1}
\hspace*{-20ex}\begin{tabbing}
\indent \= IMPA -- Estrada Dona Castorina, 110 \\
\> 22460-320 --- Rio de Janeiro --- Brazil \\
\> luis@impa.br\ \ \ \ \ - \ \ \ \  gfreitas@impa.br \end{tabbing}
}

\begin{thebibliography}{bb}  

\bibitem{alle} C. Allendoerfer,
{\it Rigidity for spaces of class greater than one}.
Amer. J. Math. {\bf 61} (1939), 633--644.

\bibitem{bbg} E. Berger, R. Bryant and P. Griffith,
{\it The Gauss equations and rigidity of isometric embeddings}.
Duke Math. J. {\bf 50} (1983), 803--892.

\bibitem{dd} M. do Carmo and M. Dajczer, 
{\it Conformal Rigidity}.
Amer. J. Math. {\bf 109} (1987), 963--985.

\bibitem{ca} E. Cartan,
{\it La d\'eformation des hypersurfaces dans l'espace euclidien
r\'eel a $n$ dimensions}.
Bull. Soc. Math. France {\bf 44} (1916), 65--99.

\bibitem{df0} M. Dajczer and L. Florit,
{\it On conformally flat submanifolds}.
Comm. Anal. Geometry {\bf 4} (1996), 261--284.

\bibitem{df1}  M. Dajczer and L. Florit,
{\it A class of austere submanifolds}.  Illinois
Math. J. {\bf 45} (2001), 735--755.

\bibitem{df2} M. Dajczer and L. Florit,
{\it Genuine deformations of submanifolds}.
Comm. Anal. Geom. {\bf 12} (2004), 1105--1129.

\bibitem{df3} M. Dajczer and L. Florit,
{\it Genuine rigidity of Euclidean submanifolds in codimension two}.
Geom. Dedicata {\bf 106} (2004), 195--210.

\bibitem{dft2} M. Dajczer, L. Florit and R. Tojeiro,
{\it On deformable hypersurfaces in space forms}.
Ann. Mat. Pura Appl. {\bf 174} (1998), 361--390.

\bibitem{dft1} M. Dajczer, L. Florit and R. Tojeiro,
{\it The vectorial Ribaucour transformation for submanifolds
and applications}.
Trans. AMS {\bf 359} (2007), 4977--4997.

\bibitem{dft3}M. Dajczer, L. Florit and R. Tojeiro,
{\it Euclidean hypersurfaces with genuine deformations
in codimension two}.
Manuscripta Math., {\bf 140} (2013), 621--643.

\bibitem{dm1} M. Dajczer and D. Gromoll,
{\it  Real Kaehler submanifolds and uniqueness of the Gauss map}.
J. Diff. Geometry {\bf 22} (1985), 13--28.

\bibitem{dm2} M. Dajczer and D. Gromoll,
{\it Isometric deformations of compact Euclidean submanifolds in
codimension two}.
Duke Math. J. {\bf 79} (1995), 605--618.

\bibitem{fera} E. V. Ferapontov,
{\it Surfaces with flat normal bundle: an explicit construction}.
Diff. Geom. Appl. {\bf 14} (2001), 15--37.

\bibitem{ft} L. Florit and R. Tojeiro,
{\it Genuine deformations of submanifolds II: the conformal case}.
Comm. Anal. Geom. {\bf 18} (2010), 397--419.

\bibitem{fz} L. Florit and W. Ziller,
{\it Manifolds with conullity at most two as graph manifolds}.
Preprint.

\bibitem{sb} V. Sbrana,
{\it Sulla variet\'a ad $n-1$ dimensioni deformabili nello spazio
euclideo ad n dimensioni}.
Rend. Circ. Mat. Palermo {\bf 27} (1909), 1--45.

\bibitem{si} S. Silva,
{\it On isometric and conformal rigidity of submanifolds}.
Pacific J. of Math. {\bf 199} (2001), 227--247.

\end{thebibliography}
\end{document}